\documentclass[11pt,twoside]{article}
\usepackage[myheadings]{fullpage}
\usepackage{amsmath,amsfonts,amssymb,bbm}
\usepackage{ulem}

\newcommand{\hide}[1]{} 
\hide{\usepackage{setspace}}
\newcommand{\nwc}{\newcommand}

\numberwithin{equation}{section}
\newtheorem{theorem}{Theorem}[section]

\newtheorem{cor}[theorem]{Corollary}

\newtheorem{lemma}[theorem]{Lemma}
\newtheorem{prop}[theorem]{Proposition}

\nwc{\todo}[1]{\smallskip
\begin{center}{TODO: #1}\end{center}
\smallskip}
\nwc{\toprove}{\qquad\mbox{YET TO PROVE}}

\nwc{\bit}{\begin{itemize}}
\nwc{\eit}{\end{itemize}}

\nwc{\qref}[1]{(\ref{#1})}
\nwc{\lbl}{\label}
\nwc{\pfrac}[2]{\left(\frac{#1}{#2}\right)}
\nwc{\pmat}[1]{ \begin{pmatrix} #1 \end{pmatrix} }
\nwc{\ip}[1]{\langle{#1}\rangle}
\nwc{\ipbig}[1]{\left\langle{#1}\right\rangle}

\nwc{\dc}[1]{#1_c}

\nwc{\bm}{\boldmath}
\nwc{\ti}{\tilde}
\nwc{\ubar}{\uline}

\nwc{\half}{{\textstyle \frac12}}
\nwc{\sfrac}[2]{{\textstyle \frac{#1}{#2}}}
\nwc{\intR}{\int_{-\infty}^\infty}
\nwc{\inv}{^{-1}}
\nwc{\D}{\partial}
\nwc{\Dt}{\partial_t}
\nwc{\Dx}{\partial_x}
\nwc{\Di}{\partial\inv}
\nwc{\R}{{\mathbb{R}}}
\nwc{\C}{{\mathbb{C}}}
\nwc{\one}{{\mathbbm{1}}}

\nwc{\rmi}{{\rm i}}
\nwc{\rme}{{\rm e}}

\nwc{\Fr}{\mathop{\rm Fr}\nolimits}
\nwc{\Bo}{\mathop{\rm Bo}\nolimits}
\nwc{\sech}{\mathop{\rm sech}\nolimits}
\nwc{\sgn}{\mathop{\rm sgn}\nolimits}
\nwc{\re}{\mathop{\rm Re}\nolimits}
\nwc{\im}{\mathop{\rm Im}\nolimits}
\nwc{\id}{\mathop{\rm id}\nolimits}

\nwc{\xpar}{\gamma}  
\nwc{\xp}{\tau}      
\nwc{\tig}{\gamma} 
\nwc{\stig}{\sqrt{\gamma}} 
\nwc{\tic}{{\tilde c}} 

\nwc{\etax}{\eta_x}
\nwc{\oneetax}{1+\eta_x^2}

\nwc{\htwonorm}[1]{\|#1\|_{H^2(\Omega)}}

\nwc{\e}{{\epsilon}}
\nwc{\ee}{{\epsilon^{1/2}}}
\nwc{\est}{{\e^{1/2}}}
\nwc{\eps}{\epsilon}
\nwc{\eph}{\epsilon} 
\nwc{\eh}{\nu}

\nwc{\eax}{e^{ax}}
\nwc{\emax}{e^{-ax}}
\nwc{\eapx}{e^{\alpha x}}
\nwc{\emapx}{e^{-\alpha x}}
\nwc{\aaw}{\alpha}
\nwc{\aw}{a}
\nwc{\lam}{\lambda}
\nwc{\lah}{\tilde\lambda}
\nwc{\gah}{\hat\gamma}
\nwc{\mesym}{{\cal M}_\eps} 
\nwc{\mosym}{{\cal M}_0} 
\nwc{\nus}{\nu} 
\nwc{\nuh}{{\hat\nu}}
\nwc{\kah}{\hat \kappa}
\nwc{\kh}{\hat k}
\nwc{\ch}{\hat c}
\nwc{\kxi}{k}

\nwc{\Zeta}{{\cal Z}}
\nwc{\Ac}{{\cal A}}
\nwc{\BB}{{\cal B}}
\nwc{\CC}{{\cal C}}
\nwc{\EE}{{\cal E}}
\nwc{\SA}{{\cal A}}
\nwc{\FF}{{\cal F}}
\nwc{\II}{{\cal I}}
\nwc{\JJ}{{\cal J}}
\nwc{\KK}{{\cal K}}
\nwc{\LL}{{\cal L}}
\nwc{\MM}{{\cal M}}
\nwc{\NN}{{\cal N}}
\nwc{\PP}{{\cal P}}
\nwc{\QQ}{{\cal Q}}
\nwc{\RR}{{\cal R}}
\nwc{\VV}{{\cal V}}
\nwc{\HH}{{\cal H}}
\nwc{\SC}{{\cal S}}
\nwc{\TT}{{\cal T}}
\nwc{\UU}{{\cal U}}
\nwc{\WW}{{\cal W}}
\nwc{\calW}{{\cal W}}
\nwc{\XX}{{\cal X}}
\nwc{\YY}{{\cal Y}}
\nwc{\GG}{{\mathcal G}}
\nwc{\Sc}{{\cal S}}
\nwc{\Dil}{\tau_{\eph}}

\nwc{\domflat}{\Omega_0}      
\nwc{\dometa}{\Omega_\eta} 
\nwc{\xf}{{\ubar x}}        
\nwc{\yf}{{\ubar y}}
\nwc{\ssf}{{\ubar s}}
\nwc{\xfrm}{{Z_1}}
\nwc{\yfrm}{{Z_2}}
\nwc{\xrm}{{z_1}}
\nwc{\yrm}{{z_2}}
\nwc{\xrmb}{{z_{10}}}
\nwc{\yrmb}{{\ubar \eta}} 
\nwc{\ueta}{{\ubar \eta}} 
\nwc{\rst}{\zeta}
\nwc{\rsti}{h} 
\nwc{\psif}{\ubar\psi}
\nwc{\Phif}{\ubar\Phi}
\nwc{\phif}{\ubar\phi}
\nwc{\zst}{{\zeta_*}}
\nwc{\hst}{{h_*}}
\nwc{\zsh}{{\zeta_{\#}}}
\nwc{\zshi}{{\zeta_{\#}\inv}}
\nwc{\hsh}{{h_{\#}}}
\nwc{\z}{{z_1}}
\nwc{\uu}{\omega}
\nwc{\us}{\theta}
\nwc{\UP}{\Upsilon}
\nwc{\uUP}{\ubar\Upsilon}
\nwc{\zf}{{\ubar z}}
\nwc{\omd}{{\varpi}}

\nwc{\Hil}{{\HH}_0}
\nwc{\Heta}{{\cal H}_\eta}
\nwc{\Hpeta}{{\cal H}'(\eta)}
\nwc{\Beta}{B(\eta)}

\nwc{\sw}{\Theta}
\nwc{\So}{{\cal S}_0}
\nwc{\Rp}{\RR_+}
\nwc{\Rm}{\RR_-}
\nwc{\AAsm}{\BB_-} 
\nwc{\AAsp}{\BB_+}
\nwc{\AAspm}{\BB_\pm} 
\nwc{\AAm}{\Ac_-}
\nwc{\AAp}{\Ac_+}
\nwc{\AApm}{\Ac_\pm}
\nwc{\BBm}{\BB_-}
\nwc{\BBp}{\BB_+}
\nwc{\AAz}{\Ac_0}

\nwc{\spy}{{C^{0,\exp}_\epsilon }}
\nwc{\spx}{{B_2^\alpha (C^{0,\exp}_\epsilon ) }}
\nwc{\spxx}{{B_3^\alpha (C^{0,\exp}_\epsilon ) }}

\nwc{\js}{{J_*}}

\nwc{\Ce}{C_\mathrm{exp}}

\nwc{\F}{{\mbox{F}}}
\nwc{\ty}{{\widetilde y}}
\nwc{\tz}{{\widetilde z}}
\nwc{\tR}{{\widetilde {\cal R}}}
\nwc{\tG}{{\widetilde  G}}

\nwc{\Lama}{\Lambda_0}
\nwc{\Lamb}{\Lambda_1}

\nwc{\rz}{|_{z=1}}
\nwc{\Dcl}{D^{\rm cl}}

\nwc{\etas}{\eta_*}
\nwc{\Us}{U_*}
\nwc{\Vs}{V_*}
\nwc{\Phis}{\Phi_*}
\nwc{\Xm}{{X_m}}
\nwc{\Ym}{{Y_m}}
\nwc{\Xmm}{{X_{m-1}}}
\nwc{\Ymm}{{Y_{m-1}}}

\nwc{\Omee}{\Omega_{\eph}}
\nwc{\Omeh}{\ti\Omega_{\eph}}
\nwc{\tiW}{\widetilde W}

\begin{document}
\title{Asymptotic linear stability of solitary water waves}
\author{Robert L. Pego\textsuperscript{1} \quad and \quad
Shu-Ming Sun\textsuperscript{2}}

\maketitle

\begin{abstract}
We prove an asymptotic stability result for 
the water wave equations linearized around small solitary waves.
The equations we consider govern irrotational flow of a fluid 
with constant density bounded below by a rigid horizontal bottom and above 
by a free surface under the influence of gravity neglecting surface tension. 
For sufficiently small amplitude waves, with waveform well-approximated
by the well-known sech-squared shape of the KdV soliton, 
solutions of the linearized equations decay at an exponential rate in an 
energy norm with exponential weight translated with the wave profile.
This holds for all solutions with no component in (i.e., symplectically orthogonal to)
the two-dimensional neutral-mode space arising from infinitesimal translational and 
wave-speed variation of solitary waves.
We also obtain spectral stability in an unweighted energy norm.
\end{abstract}

\vfil

\footnotetext[1]{Department of Mathematical Sciences and 
Center for Nonlinear Analysis, 
Carnegie Mellon University, Pittsburgh, PA 15213.
Email: rpego@cmu.edu}
\footnotetext[2]{Department of Mathematics,
Virginia Polytechnic Institute and State University,
Blacksburg, VA 24061.
Email: sun@math.vt.edu}

\pagebreak

\baselineskip=14pt
\hide{\onehalfspace}

\section{Introduction}

The discovery of solitary water waves by J.~Scott Russell in 1834 was a 
seminal event in nonlinear science.  Russell's observations gave him immediate
confidence in the significance of these waves, and led him to carry out 
an extensive program of experiments investigating solitary waves and
their interactions \cite{Russell}.
But mathematical understanding was slow to develop.
The first significant steps forward were made by 
Boussinesq \cite{Bou1871a,Bou1871b,Bou1872,Bou1877} and Rayleigh
\cite{Ray1876} by carefully balancing long-wave and small-amplitude
approximations.  
The simplest useful model (derived by Boussinesq already
in 1872, see \cite[p.~360]{Bou1877} and \cite{Miles}) 
is the famous Korteweg-de Vries equation \cite{KdV1895}.  
Its $\sech^2$ soliton solution approximates
the shape of small-amplitude solitary water waves. 

Given the status of the KdV equation as an approximate model, it is important
to understand {whether the soliton solutions of the KdV
equation are approximations of some solutions of a more exact water wave
model with similar properties}.  
In this paper, we focus on questions of stability for exact solitary wave 
solutions of the Euler equations that govern incompressible and irrotational
motions of an inviscid, constant-density fluid of finite depth.  The fluid
occupies a two-dimensional domain whose lower boundary is a flat rigid bottom
and whose upper boundary is a free surface that forms an interface with air of
negligible density and viscosity.  Surface tension on the free surface is
neglected.

For these water wave equations, the existence of solitary wave solutions with
shape well-approximated by the KdV soliton was proved by Lavrent'ev
\cite{Lav54}, Friedrichs and Hyers \cite{FrHy54} and Beale \cite{Be77}.
If the surface tension is positive and small, finite-energy, single-hump 
solitary waves are not known to exist, and indeed, 
exact traveling waves approximated by the KdV soliton may not
exist without `ripples at infinity' \cite{Be91,Sun99}.
For large surface tension, solitary water waves of depression exist \cite{Sun99},
but the relevant physical regime corresponds to water depth less 
than 0.5 cm. 

Explaining the stability of solitary water waves mathematically 
remains a very challenging problem,   
despite considerable physical and numerical evidence.
Remarkably, a valuable step forward was made
already by Boussinesq~\cite{Bou1872,Bou1877}, who argued for their stability
based on a quantity he called the `moment of instability,' which he showed
was invariant in time based on the KdV approximation.
Over a century later, Benjamin~\cite{Ben72} made use of the same quantity
as a Hamiltonian energy, constrained by a time-invariant momentum functional, 
to develop a rigorous variational method to prove orbital stability for 
the set of solitary-wave solutions of the KdV equation. 
Benjamin's arguments were improved and perfected by Bona~\cite{Bo75}.

Variational methods for orbital stability and instability in Hamiltonian wave
equations, based on the use of energy-momentum functionals, were subsequently 
greatly advanced by many authors.
Notably, the general theory of Grillakis et al.~\cite{GSS1,GSS2} 
has been applied extensively to many physical systems. 
Using variational methods of this type for the case of solitary water waves of 
depression for the Euler equations with large surface tension, orbital stability
conditional on global existence was obtained by Mielke \cite{Mi01} and Buffoni
\cite{Bu04}.  For small surface tension, such variational stability results
have also been obtained recently by Groves and Wahlen \cite{GW10} for 
oscillatory traveling wave packets of finite energy (also called solitary 
waves by several authors).

For solitary waves with zero surface tension, however, it appears hopeless to
study stability using variational methods based on constrained minimization.
As remarked by Bona and Sachs~\cite{BoSa89},
the usual energy-momentum functional is highly indefinite in this 
case---The second variation lacks the finite-dimensional indefiniteness
property key to the success of current variational methods.
Regarding the stability of solitary waves with zero surface tension, the only
existing rigorous work appears to be the recent paper of Lin \cite{Lin09}, 
which addresses the linear instability of large waves close to the wave of 
maximum height.

The present study involves a direct analysis of the Euler equations
linearized about a small-amplitude solitary wave solution.  
The linearized equations have a natural two-dimensional space of
neutral modes arising from infinitesimal shifts and changes in 
wave speed of solitary waves.  
We deduce asymptotic stability for solutions
in a space of perturbations naturally constrained to omit these neutral-mode
components, being symplectically orthogonal to them.  
Asymptotic stability is obtained in a norm that is weighted spatially to
decay exponentially behind the wave profile.  The time decay of such a
norm corresponds to unidirectional scattering behavior for wave perturbations.
The weighted-norm linear stability analysis is also used to obtain a spectral
stability result in an unweighted energy norm.
Our main results are stated precisely in section 3.  

The use of exponential weights to obtain nonlinear asymptotic orbital
stability for solitary waves was developed for KdV solitons by Pego and
Weinstein~\cite{PW94}, for regularized long-wave equations by Miller and
Weinstein~\cite{MiWe96} and for Fermi-Pasta-Ulam lattice equations by Friesecke
and Pego~\cite{FP1,FP2,FP3,FP4}.  Finiteness of an exponentially weighted norm
imposes a condition of rapid decay in front of the wave profile. But Mizumachi
\cite{MP08,MZ09} recently showed how to prove asymptotic orbital stability for
FPU solitary waves perturbed in the energy space, by using exponential weights
together with dispersive wave propagation estimates as developed by Martel and
Merle.

Nonlinear stability for solitary water waves remains an open problem.
This issue would likely involve a general
global existence theory for small-amplitude 2D fluid motions,
which is not yet available despite the substantial progress 
on well-posedness questions by Wu \cite{Wu2d,Wu09}. 

There are a number of other works on (in)stability for 2D solitary water
waves that concern the case of waves of depression with large surface tension.
These include results on 2D spectral stability for finite-wavelength
perturbations \cite{HS02a}, spectral instability for transverse (3D)
perturbations \cite{PSun}, and a full nonlinear instability result
for 3D perturbations by Rousset and Tzvetkov \cite{RT}.

A convenient tool for singular perturbation theory, 
used in \cite{PSun} and in the
present paper to study spectrum in the KdV scaling limit of long time and length
scales, is an operator-theoretic generalization of Rouch\'e's theorem due to
Gohberg and Sigal \cite{GS71}.  This use of the KdV scaling contrasts
with works by Craig \cite{Craig} and Schneider and Wayne \cite{ScWa00b} that
concern the validity of the KdV approximation for water waves over time scales
of order $O(\eph^{-3})$ for waves of amplitude $O(\eps^2)$ that are long with
length scales of order $\eph\inv$.  Our use of the KdV approximation occurs in
the spectral domain, where it is used to obtain partial information regarding
the behavior of solutions to the linearized equations in the limit
$t\to\infty$.
To establish stability for time and space scales unrelated to the regime of
validity of the KdV approximation requires a different technique for 
dealing with the linearized Euler equations, which resemble a wave equation
with variable coefficients.  We develop a method that obtains
resolvent bounds from symmetrized weighted-norm energy estimates that use
Fourier filters to cut off low frequencies.

\vfil

\section{Equations of motion and eigenvalue problem}

In this section, we derive the equations of motion linearized around a
solution steady in a frame 
moving at a constant speed $c$ to the right, 
and formulate the associated eigenvalue problem.

{\bf Basic equations.} 
We deal with an inviscid, incompressible and irrotational fluid of constant
density $\rho$ that is bounded above by a free surface $y = \eta (x,t)$
and below by a horizontal rigid bottom $y =-h$. 
The velocity field $(u,v)$ is related to the velocity potential $\phi$
and the stream function $\psi$ by
\begin{equation}\lbl{D.uv}
(u,v)=(\phi_x,\phi_y)=(\psi_y,-\psi_x) .
\end{equation}
On the free surface $y=\eta(x,t)$, the kinematic and Bernoulli
equations are:
\begin{eqnarray}
\Dt\eta + u \eta_x &=& v, 
\label{e.eta1}\\
\Dt\phi + \half(\phi_x^2+\phi_y^2) + g\eta &=&0.
\label{e.phi1}
\end{eqnarray}
To make the problem non-dimensional, we let
\begin{equation}
(x, y, t) =(h\ti x , h\ti y, h\ti t/c), 
\qquad (\eta,u,v,\phi,\psi) = (h\ti\eta, c\ti u, c\ti v,ch\ti\phi,ch\ti\psi).
\end{equation}
After dropping the tildes, 
the equations take again the same form in the non-dimensional variables,
with $g$ replaced by 
\begin{equation}
\tig=\frac{gh}{c^2}= \frac1{\Fr^2}, 
\end{equation}
where $\Fr=c/\sqrt{gh}$ is the Froude number.

In the fluid region, where now $-1<y<\eta(x,t)$, $-\infty<x<\infty$, the velocity
potential and stream function are harmonic and are taken to satisfy the no-penetration boundary
conditions
\begin{equation}\lbl{E.bottom1}
\phi_y(x,-1) = 0, \quad  \psi(x,-1) = 0 \quad (-\infty<x<\infty).
\end{equation}
The dynamics is described in terms of the surface traces defined by 
\begin{equation}\lbl{D.Phi}
\Phi(x,t)=\phi(x,\eta(x,t),t), \quad
\Psi(x,t)=\psi(x,\eta(x,t),t).
\end{equation}
Then
\begin{equation}\lbl{D.UV}
\pmat{U\cr V}:=\pmat{ \Phi_x\cr\Psi_x} = 
\pmat{ \phi_x+\etax\phi_y \cr \psi_x+\etax\psi_y}
= 
\pmat{ 1&\etax\cr\etax&-1}
\pmat{u\cr v},
\end{equation}
and we will write 
\begin{equation}\lbl{D.Mmat}
\qquad
M(\etax)= 
\pmat{ 1&\etax\cr\etax&-1},
\qquad  
M(\etax)\inv
=\frac{M(\etax)}{1+\etax^2}.
\end{equation}

The non-dimensional equations of motion now take the form
\begin{eqnarray}
\Dt\eta &=& v-\etax u = -V, 
\label{e.eta2} \\
\Dt\Phi &=& 
\Dt\phi+\phi_y\Dt\eta  = -\tig\eta-\half(u^2+v^2)-v(v-\etax u)
\nonumber\\ 
&=&
-\tig\eta - \half (U,V)M(\etax)\inv(U,V)^T.
\label{e.phi2} 
\end{eqnarray}
After transforming to a moving frame with $\hat x=x-t$ ($=(x-ct)/h$
dimensionally) and dropping the hats, the time derivative $\Dt$ is 
replaced by $\Dt-\Dx$. 
A solitary wave is a steady solution of the resulting equations.

It is convenient to regard the wave motion as determined by the
evolution of the pair $(\eta,\Phi)$, with $\Psi$ and $V=\Psi_x$ determined 
from $(\eta,\Phi)$ by solving for the stream function using
Laplace's equation and the relevant boundary conditions, namely
(suppressing the $t$ variable)
\begin{eqnarray}
&&\psi_{xx}+\psi_{yy}=0 \qquad (-\infty<x<\infty,\ -1<y<\eta(x)), 
\lbl{E.psi1}\\
&&\psi(x,-1)=0, \quad \psi_y-\etax\psi_x = U(x) 
\qquad (-\infty<x<\infty,\ y=\eta(x)).
\lbl{E.psi2}
\end{eqnarray}
We write 
\begin{equation}\lbl{D.BH}
\Psi=\Heta\Phi=\psi(x,\eta(x)), \qquad V= \Psi_x.
\end{equation}
Up to a normalization, $\Heta$ is a Hilbert transform 
for the fluid domain. (Note $\phi+i\psi$ is an analytic function
of $x+iy$.)
This map will be studied in detail in a later section.

{\bf Linearization.}
We linearize the equations in the moving frame about a steady solution, 
denoting linearized variables with a dot.
These linearized equations of motion take the form
\begin{eqnarray}
\lbl{E.lin1}
0&=& (\Dt-\Dx) \dot\eta + \dot V  , \\
0&=& (\Dt-\Dx) \dot\Phi + \tig \dot\eta + u\dot U + v\dot V
-uv\Dx\dot\eta.
\lbl{E.lin2}
\end{eqnarray}
Of course $\dot U=\Dx\dot\Phi$. To relate $\dot V$ to $(\dot\eta,\dot \Phi)$,
we linearize the boundary-value problem \qref{E.psi1}-\qref{E.psi2}
by formally differentiating with respect to a variational parameter.
The variation $\dot\psi$ is harmonic in the fluid domain, zero on the bottom,
and on the free surface $y=\eta(x)$ satisfies
\[
\dot U = 
\dot\psi_y-\etax\dot\psi_x
-\dot\eta_x\psi_x+(\psi_{yy}-\eta_x\psi_{xy})\dot\eta.
\]
Since $\psi_{yy}-\eta_x\psi_{xy}=-\Dx( \psi_x(x,\eta(x)))$ and
$-\psi_x=v$ , this means
\begin{equation}\lbl{E.dpsi2}
\dot\psi_y-\etax\dot\psi_x = \Dx(\dot \Phi(x) -\dot \eta(x) v(x,\eta(x))).
\end{equation}
and by \qref{D.BH} this means 
$\dot\psi(x,\eta(x)) = \Heta(\dot\Phi-v\dot\eta)$.
Hence $\dot V=\Dx\dot\Psi$ where
\begin{equation}\lbl{E.dV}
\dot \Psi = \psi(x,\eta(x))\,\dot{} =  
  \dot\psi(x,\eta(x))+\psi_y(x,\eta(x))\dot\eta 
= \Heta (\dot\Phi-v\dot\eta ) +u\dot\eta.
\end{equation}

We have found it to be important (much more than merely convenient) 
to study the linearized equations of motion in terms of the 
combination of $\dot\eta$ and $\dot\Phi$ expressed as
\begin{equation}\lbl{D.dphi}
\dot\phi = \dot\Phi-v\dot\eta.
\end{equation}
This is the 
{\it surface trace of the variation of velocity potential},
rather than the variation of the surface trace. 
A similar observation was made by Lannes \cite{La05} in his
treatment of well-posedness for 3D water waves locally in time.
In terms of the pair $(\dot\eta,\dot\phi)$, the linearized 
equations of motion take the form
\begin{equation}
\lbl{e.lin2}
(\Dt-\Ac_\eta)
\pmat{\dot\eta \cr \dot\phi} = 0,
\qquad
\Ac_\eta = \pmat{ \Dx(1-u) & -\Dx\Heta \cr  
-\tig+(1-u)v'  & (1-u)\Dx },
\end{equation}
where $v'$ is the multiplier $v'(x)=\Dx(v(x,\eta(x)))$.
Our analysis will show that the initial-value
problem for the linear system \qref{e.lin2}
is well-posed and (conditionally) asymptotically stable in a certain 
weighted function space. The components $\dot\eta$ and $\dot\phi$ 
will belong to spaces of different order, however, and 
this complicates the problem of studying stability questions
directly using the variables $(\dot\eta,\dot\Phi)$.

{\bf Eigenvalue problem.}
Looking for solutions of \qref{e.lin2} in the form
$(\dot\eta,\dot\phi)=e^{\lambda t}(\eta_1(x),\phi_1(x))$
leads to the associated eigenvalue problem
\begin{equation}
\lbl{E.evp1}
\pmat{ \lam-\Dx(1-u) & \Dx\Heta \cr  
\tig-(1-u)v'  & \lam-(1-u)\Dx }
\pmat{\eta_1 \cr \phi_1} = 0.
\end{equation}
The hardest part of our analysis of the linearized dynamics 
involves showing that, in an appropriate function space,
this equation has no nontrivial solutions 
for all nonzero $\lam$ in a half plane $\re\lam\ge -\beta$
for some $\beta>0$ depending on the wave amplitude.

\section{Main results}

Our main result is an asymptotic linear stability result 
for the classical family of small-amplitude solitary water waves 
that exist for Froude number slightly more than 1, meaning $\tig<1$.
Asymptotic stability is conditional on the absence of neutral-mode 
components arising from translational shifts of the solitary wave,
and wave-speed variation, as is standard.
The precise results involve $L^2$ spaces with
exponential weights $e^{ax}$ that decay to the left (having $a>0$).
For $a\in\R$, we define $L^2_a$ to be the Hilbert space 
\[
L^2_a=\{f\mid e^{ax}f\in L^2(\R)\},
\]
with inner product and norm
\[
\ip{f,g}_a = \intR f(x)\overline{g(x)} e^{2ax}\,dx,
\qquad \|f\|_a = \| e^{ax}f\|_{L^2}.
\]
Also, $H^s_a = \{f\mid e^{ax}f\in H^s(\R)\}$ 
will denote a weighted Sobolev space with norm 
that is expressed in terms of the Fourier transform 
$\FF f(k)=\hat f(k)=\intR e^{-ikx}f(x)\,dx$ as
\[
\|f\|_{H^s_a}=\|\eax f\|_{H^s} = 
\left(\frac1{2\pi}\intR (1+k^2)^s |\hat f(k+ia)|^2\,dk\right)^{1/2}.
\]

{\bf Group velocity and weighted norms.} 
The use of these weighted norms is motivated as follows.
For linearization at the trivial solution $\eta=\Phi=0$,
the Hilbert transform for the fluid domain is the Fourier multiplier
$\HH_0=i\tanh D$ with $L^2$ symbol $i\tanh k$ (see section~\ref{S.hil}).
Then the dispersion relation for solutions of \qref{e.lin2}
with space-time dependence $e^{ikx-i\omd t}$ is 
\begin{equation}\lbl{e.disp0}
\omd = -k \pm \sqrt{\tig k\tanh k}.
\end{equation}
The solitary waves that we study travel faster than the speed of long
gravity waves, meaning $c>\sqrt{gh}$ and so $\tig<1$.
Thus, in this regime the group velocity 
of linear waves (relative to the solitary wave) is always negative:
\begin{equation}
\frac{d\omd}{d k }<0.
\end{equation}
Heuristically, linear waves scatter to the left. Our analysis
makes essential use of this directionality 
by measuring perturbation size using weights $e^{ax}$ with $a>0$.
As a simple example, the solution $\eta(x,t)=f(x+t)$ 
of the transport equation $\Dt\eta=\Dx\eta$ 
satisfies $\|\eta(\cdot,t)\|_a= e^{-t}\|f\|_a$. 

In analytic terms, the isomorphism $f\mapsto \eax f$ from $L^2_a$ to $L^2$
maps a Fourier multiplier $\Ac(D)$ acting on $L^2_a$ 
to the weight-transformed operator $\eax\Ac(D)\emax=\Ac(D+ia)$
acting on $L^2$. 
The $L^2$-symbol $\Ac(k)$ of the former is shifted to 
the symbol $\Ac(k+ia)$ of the latter.
The $L^2_a$-spectrum of $\Ac(D)$ is the 
closure of the image of the latter symbol. 
This is so since the resolvent $(\lam-\Ac(D))\inv$ is bounded in $L^2_a$
exactly when the map $f_a\mapsto (\lam-\Ac(k+ia))\inv \hat f_a(k)$
is bounded in $L^2$, where $f_a=\eax f$.
For the Fourier multipliers 
\[
\Ac_\pm(D)= iD\pm \sqrt{-\tig D\tanh D},
\]
which correspond to the branches of the dispersion relation
\qref{e.disp0} for our water-wave problem,
the $L^2_a$-spectrum shifts from the
imaginary axis into the left half-plane for small $a>0$ exactly because the 
relative group velocity is negative.
The same idea underlies the use of weights to obtain nonlinear
asymptotic stability for solitary waves of the KdV equation in
\cite{PW94} and of FPU lattice equations in \cite{FP2,FP3,FP4}.

{\bf Energy and weighted norms.} 
Zakharov~\cite{Za68} showed that the water wave equations
have a canonical Hamiltonian structure in terms of $(\eta,\Phi)$ with
(nondimensional) Hamiltonian
\begin{equation}\lbl{e.ham}
\frac12 \intR\int_{-1}^{\eta(x)} |\nabla\phi|^2\,dy\,dx+ 
\frac12\intR \tig\eta^2 \,dx 
=\frac12\intR \left(\Phi(-\Dx\Heta)\Phi+\tig\eta^2\right) dx.
\end{equation}
The space that we use to study asymptotic stability of the linearized system
\qref{e.lin2} is equivalent to a weighted linearization of this
Hamiltonian about a flat surface. Namely, stability will be studied
with $(\dot\eta,\dot\phi)$ 
in the space $Z_a=L^2_a\times H^{1/2}_a$ with norm equivalent to the norm of
$(\dot\eta,\sqrt{D\tanh D}\dot\phi)$ in $L^2_a\times L^2_a$.

{\bf Scaling.} We study waves in the regime where the parameter
\begin{equation}\lbl{d.eph}
\eph = \sqrt{1-\tig}
\end{equation}
is small and positive. For all $\eph$ in this well-studied regime, 
there is an even solitary-wave surface elevation $\eta$ with $\eta$ and 
surface velocity $(u,v)$ approximately given by 
\[
\eta(x)\sim u(x)\sim \eph^2\sw(\eph x), \qquad
v(x)\sim -\eph^3\sw'(\eph x)
\]
where
\begin{equation}\lbl{d.S0}
\sw(x)= \sech^2(\sqrt3 x/2).
\end{equation}
For precise statements with estimates we use, see
Theorems~\ref{t.sw} and \ref{t.swa}.
The significance of these results is that we 
use stability information for the KdV soliton with the profile \qref{d.S0}
to study the eigenvalue problem \qref{E.evp1} for $|\lam|$ small,
using the KdV scaling $\hat x=\eph x$, $\lam=\eph^3\lah$. 
Because of this scaling, we take the weighted-norm exponent
to have the form $\aw=\eph\aaw$, where $\aaw$ is 
required to satisfy $0<\aaw<\sqrt3$ to have $\sw\in L^2_\aaw$. 
For convenience in analysis, our stability results are formulated with the
tighter restriction $0<\aaw\le\frac12$.  
The parameter $\aaw$ is taken as any fixed number in this range.

{\bf Neutral modes.} 
The solitary waves we study belong to a two-parameter family, 
smoothly parameterized by translation and Froude number (equivalently
translation and wave speed c). 
By consequence, as usual the value $\lambda=0$ is an eigenvalue of
$\Ac_\eta$ with algebraic multiplicity two, 
with generalized eigenfunctions produced by 
differentiation with respect to $x$ and $c$.  Denoting these
functions with the notation
\[
z_x=\pmat{\eta_x\cr \phi_x}, \qquad z_c = \pmat{\eta_c\cr \phi^+_c},
\]
we have $\Ac_\eta z_x=0$, $-\Ac_\eta z_c=z_x$.  The details are developed
in Appendix B. (The notation $\phi^\pm_c$ indicates that different choices
of an integration constant are made to ensure 
$\Phi^\pm=\D_x\inv U\in H^{1/2}_{\pm a}$.)

Solutions of \qref{e.lin2} that lie in the neutral-mode space
spanned by $z_x$ and $z_c$ do not decay in time, naturally.
A necessary condition that a solutions of \qref{e.lin2} decay in time
is that it should have no component in this neutral-mode space.
The precise spectral meaning of this 
(being annihilated by the spectral projection for the eigenvalue $\lam=0$)
can be expressed in a simple form, 
due to the canonical Hamiltonian structure of the problem.
Namely, it turns out to be necessary that the solution be
{\it symplectically orthogonal}
to the neutral mode space, meaning that 
\begin{equation}\lbl{c.symp}
0 = \intR \dot\eta \phi_x - \dot\phi \eta_x \,dx,
\qquad
0 = \intR \dot\eta \phi^-_c - \dot\phi \eta_c \,dx.
\end{equation}

{\bf Results.}
Our main results concern asymptotic stability for the linearized equations 
in a weighted norm, and spectral stability in an unweighted norm.

\begin{theorem}\lbl{t.main}
(Asymptotic stability with weights)
Fix $\aaw\in(0,\frac12]$ and set $a=\aaw\eph$. 
If $\eph>0$ is sufficiently small and
$\eta$, $u$, $v$ correspond to the solitary wave profile given by
Theorem~\ref{t.sw}, then the following hold.
\begin{itemize}
\item[(i)]
With domain $H^1_a\times H^{3/2}_a$, $\Ac_\eta$
is the generator of a
$C^0$-semigroup in $Z_a=L^2_a\times H^{1/2}_a$.
\item[(ii)] Whenever $\re\lam\ge-\frac16\aaw\eph^3$ and $\lam\ne0$, 
$\lam$ is in the resolvent set of $\Ac_\eta$.
\item[(iii)] The value $\lam=0$ is a discrete eigenvalue of $\Ac_\eta$
with algebraic multiplicity 2.
\item[(iv)] There exist 
constants $K>0$ and $\beta>\frac16\aaw\eph^3$ depending on $\eps$ and
$\aaw$, such that for all $t\ge0$,
\[
\|\exp(t\Ac_\eta)\dot z\|_{Z_a} \le K e^{-\beta t}\|\dot z\|_{Z_a},
\]
for every initial state $\dot z=(\dot\eta,\dot\phi)$
that satisfies the symplectic orthogonality conditions \qref{c.symp}.
\end{itemize}
\end{theorem}

\begin{theorem}\lbl{t.nowt} 
(Spectral stability without weights)
For $\eph>0$ sufficiently small, in the space of pairs
$(\eta_1,\phi_1)$ such that 
\begin{equation}\lbl{e.espc}
\intR 
\phi_1(D\tanh D)\phi_1 + \eta_1^2 \, dx <\infty
\end{equation}
the spectrum of the operator $\Ac_\eta$ 
is precisely the imaginary axis.
\end{theorem}

The asymptotic stability statement in part (iv) of Theorem~\ref{t.main} will
be proved as a consequence of the Gearhart-Pr\"uss spectral mapping theorem
\cite{Pr85} by establishing that the operator $\Ac_\eta$ has uniformly
bounded resolvent $(\lam-\Ac_\eta)\inv$ for $|\lam|$ large with
$\re\lam\ge-\frac16\aaw\eph^3$, and using parts (ii) and (iii) to infer
the resolvent restricted to the spectral complement of the generalized
kernel is uniformly bounded for all
$\lam$ with $\re\lam\ge-\frac16\aaw\eph^3.$
The spectral stability result in Theorem~\ref{t.nowt} 
is proved using Theorem~\ref{t.main} and 
symmetry of the problem without weights under space and time reversal.


\section{Riemann mapping and Hilbert transform for the fluid domain}
\lbl{S.hil}

\subsection{Riemann stretch and strain}
We will make much use of a Riemann mapping from the 
fluid domain $\dometa$ to the flat strip $\domflat$, with
\begin{align*}
&\dometa = \{(x,y): -\infty<x<\infty, \ -1<y<\eta(x)\},
\\
&\domflat = \{(\xf,\yf) : -\infty<\xf<\infty, \ -1<\yf<0\}.
\end{align*}
To denote the corresponding Riemann mapping and its inverse, we write
\begin{equation} \label{rm.1}
(\xf,\yf)=(\xfrm(x,y),\yfrm(x,y)), 
\qquad
(x,y) = (\xrm(\xf,\yf),\yrm(\xf,\yf)).
\end{equation}
A key quantity is the `Riemann stretch' $\rst$ at the surface, 
given (with its inverse $\rsti=\rst\inv$) by 
\begin{equation}
 \rst(\xf):= \xrm(\xf,0),
\qquad 
 \rsti(x)= \xfrm(x,\eta(x)). 
\end{equation}
The function $\yrm$ is harmonic in the strip $\domflat$, with boundary 
conditions
\begin{equation}
\yrm(\xf,-1)=-1, \qquad \yrm(\xf,0)=\yrmb(\xf):=\eta\circ\rst(\xf).
\label{z2bc}
\end{equation}
Taking the Fourier transform in $\xf$ leads to the formula
\begin{equation}
\yrm(\xf,\yf) = \yf + 
\frac{1}{2\pi}\intR e^{ik\xf}
\frac{\sinh k(\yf+1)}{\sinh k}
\intR e^{-iks} 
\yrmb(s) \,ds\,dk.
\end{equation}
Using the Cauchy-Riemann equation $\D_\xf\xrm=\D_\yf\yrm$,
we find that the `Riemann strain' defined by $\omega=\zeta'-1$
satisfies
\begin{equation}\lbl{d.omega}
\omega(\xf) = \rst'(\xf)-1 = D\coth D\, \yrmb(\xf).
\end{equation}
Integrating in $\xf$ with an arbitrary constant of integration,
we find that we can write
\begin{equation}\label{rm.z10}
\rst(\xf) = \xf-i\coth D\,\yrmb(\xf)+ c_0= \xf - \int_\xf^\infty\yrmb(s)\,ds
+ Q_1(D) \yrmb(\xf) +c_0,
\end{equation}
where $Q_1(D)$ is the Fourier multiplier with symbol bounded on $\R$ 
given by
\[
Q_1(k) = \frac{k\cosh k - \sinh k}{ik\sinh k} = i(k\inv-\coth k).
\]
If $\eta$ is given, then since $\ueta=\eta\circ\zeta$,
Eq.~\qref{rm.z10} is a fixed-point equation that
determines $\rst$ and therefore $\rsti=\rst\inv$.
It will turn out to be more convenient in our analysis, however, 
to directly study the Riemann strain $\omega$,
and recover other quantities such as $\zeta$ and $\eta$ from this.

\subsection{Hilbert transform} 
The operator $\Heta$ admits a convenient expression 
in terms of the Riemann stretch $\zeta$. 
To see this, first we introduce pullback operators 
$\zsh$ and $\zst$ via
\begin{equation}\label{D.pull}
\zsh U(\xf)= U\circ\zeta(\xf), \qquad
\zst U(\xf) =  \zeta'\zsh U(\xf) = 
(U\circ\zeta)(\xf)\, \zeta'(\xf) .
\end{equation}
For later use, note that since $h=\zeta\inv$, 
the chain rule yields the simple relations
\begin{equation}\label{E.zrel}
\D\zsh=\zst\D, \quad
\zst = \zeta'\zsh= \zsh(1/h'), \quad
\hst\zst=\zst\hst=\id.
\end{equation}
Write $\psif(\xf,\yf)=\psi(x,y)$.
Then $\psif$ is harmonic in $\domflat$, and the boundary conditions
\qref{E.psi2} transform to 
\begin{equation}\label{psifbc}
\psif(\xf,-1)=0, \qquad
\D_\yf\psif(\xf,0)=  \zst U(\xf)=\D_\xf(\Phi\circ\zeta)(\xf).
\end{equation}
By Fourier transform we find
\[
\psif(\xf,\yf) = 
\frac{1}{2\pi}\intR e^{ik\xf}
\frac{\sinh k(\yf+1)}{k\cosh k}
\intR e^{-iks} \D(\Phi\circ\zeta)(s)\,ds\,dk,
\]
so since $\psi(x,\eta(x))=\psif(h(x),0)$, after
an integration by parts we find
\begin{equation}\label{e.B1}
\Psi(x)=\psi(x,\eta(x)) = 
\frac{1}{2\pi}\intR e^{ik\rsti(x)}(i\tanh k)
\intR e^{-iks}\Phi\circ\zeta(s)\,ds\,dk.
\end{equation}
In other words, since $\Psi=\Heta \Phi$ we have 
(with $\FF$ denoting Fourier transform)
\begin{equation}\label{e.B2}
\Heta  = \zshi \Hil \zsh ,
\qquad \Hil = i\tanh D =  \FF\inv (i\tanh k) \FF.
\end{equation}
Here $\Hil$ is the Hilbert transform for the top boundary of the 
strip. Though we will make no use of the fact, it is explicitly given 
in terms of an integral kernel by 
\begin{equation}\label{e.hil}
\Hil U(x) = \intR k_0(x-s)U(s)\,ds, \qquad
k_0(x) = \frac{-1} {2\sinh(\pi x/2)}.
\end{equation}

\subsection{Linearization}
To justify the later calculation of generalized eigenmodes (in Appendix B), we 
explain here how the formal linearization formula in \qref{E.dV}
follows from the representation formula in \qref{e.B2}. 

To proceed, start with a family of (smooth) Riemann strains
$\omega$ small in $L^2\cap L^2_a$ and depending smoothly on a
variational parameter, and compute 
\[
\zeta(\xf)=\xf+\D\inv \omega + c_0, 
\quad\ueta=\frac{\tanh D}{D}\omega,
\quad  \eta=\ueta\circ \zeta\inv.
\]
Determine conjugate harmonic functions $z_1$, $z_2$ in the 
strip $\domflat$ such that \qref{z2bc} holds and $z_1(\xf,0)=\zeta(\xf)$.
Then the function $\Zeta(\xf+i\yf)=z_1(\xf,\yf)+iz_2(\xf,\yf)$
yields the Riemann mapping 
of $\domflat$ to $\dometa$ as described above. 

Also take a family of (smooth) functions $\Phi$ 
(free surface velocity potential)
and introduce $\phif$ as the harmonic extension of $\Phi\circ\zeta$
into $\domflat$ satisfying $\D_\yf\phif=0$ at $y=-1$, and $\psif$ as the 
harmonic function conjugate to $\phif$ and satisfying $\psif=0$ at $y=-1$.
Then
\[
\phif(\xf,\yf)= \frac1{2\pi}
\intR e^{ik\xf} \frac{\cosh k(\yf+1)}{\cosh k} 
\intR e^{-iks} \Phi\circ\zeta(s)\,ds\,dk.
\]
Write 
\[
\uUP(\xf+i\yf)=\phif(\xf,\yf)+i\psif(\xf,\yf),
\qquad \UP=\uUP\circ\Zeta\inv. 
\]
Regarding $\domflat$ as a subset of the complex plane, 
we have that $\uUP$ and $\Zeta$ are analytic in $\domflat$ 
and that $\UP$ is analytic in $\dometa=\Zeta(\domflat)$. 
Define $\phi$ and $\psi$ to satisfy
\[
\UP(x+iy)=\phi(x,y)+i\psi(x,y), \qquad (x,y)\in\dometa.
\]
$\Psi(x)=\psi(x,\eta(x))$ is the trace on the fluid surface and satisfies 
\[
\Psi\circ\zeta = i\tanh D\, (\Phi\circ\zeta),
\]
due to $\Psi\circ\zeta(\xf)=\psif(\xf,0)$ and 
$\Phi\circ\zeta(\xf)=\phif(\xf,0)$
and the boundary condition $\im \uUP=0$ at $y=-1$.

Denoting the derivative with respect to the variational parameter by a dot, 
we have 
\begin{equation}\lbl{e.dot1}
\dot\Psi\circ\zeta+ \dot\zeta\Psi_x\circ\zeta  = 
i\tanh D
(\dot\Phi\circ\zeta+ \dot\zeta\Phi_x\circ\zeta ) 
\end{equation}
Note that $\dot\Zeta \UP'\circ \Zeta $ is analytic in $\domflat$ and 
has zero imaginary part on the bottom $y=-1$.
This means that the real and imaginary parts of the surface trace are related 
by the Hilbert transform for the strip. But 
$\dot\Zeta=\dot z_1+i\dot z_2= \dot\zeta+i\dot\ueta$ on $y=0$, whence
(abusing notation to write $\psi_x$ for $\psi_x(x,\eta(x))$ with 
$x=\zeta(\xf)$, etc.)
\begin{equation}\lbl{e.dot2}
\psi_x \dot\zeta+\psi_y\dot\ueta = (i\tanh D)(\phi_x\dot\zeta+\phi_y\dot\ueta).
\end{equation}
Using the formulas 
\[
\dot\ueta=\dot\eta\circ\zeta+\dot\zeta\eta_x\circ\zeta,
\qquad \Phi_x = \phi_x+\phi_y\eta_x,
\qquad \Psi_x = \psi_x+\psi_y\eta_x,
\]
together with \qref{e.dot1} now yields
\[
\dot\zeta \Psi_x = -\dot\eta\psi_y+ \dot\zeta\psi_x+\dot\ueta\psi_y
\]
and similarly for $\Phi$. Combining this with \qref{e.dot2} yields
\begin{equation}
\dot\Psi\circ\zeta - \dot\eta \psi_y = 
i\tanh D (
\dot\Phi\circ\zeta - \dot\eta \phi_y),
\end{equation}
and composing with $h=\zeta\inv$ yields the desired linearization formula \qref{E.dV}:
\begin{equation}
\dot\Psi - u\dot\eta = 
\Heta (
\dot\Phi - v\dot\eta) 
\end{equation}
(with $(u,v)=(\phi_x,\phi_y)=(\psi_y,-\psi_x)$ on $y=\eta(x)$).

\section{Solitary wave profiles}

In this section and Appendix A, 
we will give a simple self-contained account of the existence
of small solitary waves by fundamentally the same approach as
Friedrichs and Hyers~\cite{FrHy54}, establishing the estimates 
that we need regarding convergence of the scaled wave profiles
in the KdV limit. 

First, note that from \qref{e.eta2}-\qref{e.phi2},
the steady equations for a solitary wave are 
\begin{equation}
\D_x\eta = V =\D_x\Heta \Phi, \qquad 
U-\gamma\eta = \frac12 (U,V)M(V)\inv(U,V)^T 
= \frac{U^2-V^2+2UV^2}{2(1+V^2)},  
\label{e.Us1}
\end{equation}
whence
\begin{eqnarray}
&& \eta = \Heta \Phi,
\qquad 
U-\gamma\eta = \frac12(U^2-V^2)+\gamma\eta V^2. 
\label{e.Us}
\end{eqnarray}
Using \qref{e.B2} and changing
variables by applying $\zsh$ we must have 
$\ueta = i\tanh D(\Phi\circ\zeta)$, hence
by \qref{d.omega},
\begin{equation}\label{e.bareta}
\uu = \zeta'-1 = (D\coth D)\ueta =  
\D (\Phi\circ\zeta) = \zeta' U\circ\zeta.
\end{equation}
Then we find
\begin{equation}\lbl{e.zshUV}
U\circ\zeta = \frac{\uu}{1+\uu} = \uu - \frac{\uu^2}{1+\uu},
\qquad 
V\circ\zeta = \frac{\D\ueta}{\zeta'}= \frac{i\tanh D\, \uu}{1+\uu}. 
\end{equation}

It is convenient to apply $\zsh$ to (\ref{e.Us}b), and isolate
$\uu$ on the left-hand side. This turns (\ref{e.Us}b) into a
fixed-point equation for the Riemann strain $\uu$, in the form
\begin{equation}\label{e.uufp}
\uu = \left(1-\gamma\frac{\tanh D}D\right)\inv
\left( 
\frac{\frac32\uu^2 + \uu^3 - \frac12 (i\tanh D\,\uu)^2(1-2\gamma \yrmb ) }
{(1+\uu)^2}
\right), \qquad
\ueta = \frac{\tanh D}{D}\uu.
\end{equation}
The following result provides scaled bounds
for the fixed point approximating the 
$\sech^2$ KdV profile from \qref{d.S0}.
The proof is given in appendix A.

\begin{theorem}\lbl{t.sw}
Let $\aaw\in(0,\sqrt3)$, $m\ge2$, $\nus\in(0,1)$,
and $\sw(x)=\sech^2(\sqrt3 x/2)$.
Then for $\eph>0$ sufficiently small, equation~\qref{e.uufp}
has a unique even solution in $H^1_a$ of the form 
\begin{equation}\lbl{d.us}
\uu(\xf)=\eph^2\us(\eph\xf)
\end{equation}
with $\|\us-\sw\|_{H^m_\aaw} <\eph^\nus$.
Moreover, the map $\eph\mapsto \uu$ is smooth.
\end{theorem}

The coefficients that appear in the linearized
system \qref{E.lin2} can now be expressed as follows.
Using \qref{D.UV}, \qref{e.Us1} and \qref{e.zshUV}, 
on the fluid surface we have the formulas
\begin{equation}\lbl{e.uva}
u = \frac{U+V^2}{1+V^2}, \qquad
v = \frac{-V(1-U)}{1+V^2},
\end{equation}
\begin{equation}\lbl{e.uvomeg}
\zsh u = \frac{ \uu+\uu^2+ (i\tanh D\uu)^2}{(1+\uu)^2+(i\tanh D\uu)^2},
\quad
\zsh v = \frac{ -i\tanh D\uu}{(1+\uu)^2+(i\tanh D\uu)^2},
\quad
\zsh v'=\frac{\D\zsh v}{1+\uu}.
\end{equation}

Note that $\ueta$, $\zsh u$ and $\zsh v'$ are even functions
(since functions have the same parity as their Fourier transform).
With the choice $c_0=\int_0^\infty\ueta(s)\,ds$, $\zeta$ is odd,
and $\eta$ and $u$ are even, with $v=(u-1)\eta_x$ odd.
Formally, we have the leading order approximations
\begin{equation}
\uu(\xf) \sim \ueta(\xf) \sim \zsh u(\xf)\sim \eph^2 \sw(\eph \xf), 
\qquad
\zsh v'(\xf) \sim -\eph^4 \sw''(\eph \xf) 
\end{equation}
For making estimates involving the quantities in \qref{e.uvomeg} it is 
useful to note that
\begin{equation}\lbl{i.taus}
\|i\tanh\eph D\us\|_{H^1} = 
\left\| \eph \D \frac{\tanh\eph D}{\eph D}\us\right\|_{H^1}
\le \eph \|\us\|_{H^2}.
\end{equation}

\section{Transforming the system}

{\bf Flattening.} 
Given the form of $\Heta$ in \qref{e.B2}, 
it appears convenient to transform the eigenvalue problem in
\qref{E.evp1} to work in variables associated with the 
flattened domain.  We make a similarity transform of \qref{E.evp1} by 
applying the operator $\zst=\zeta'\zsh$ from \qref{D.pull}
to the first equation and $\zsh$ to the second,
introducing the variables
\begin{equation}\lbl{d.eU2}
\pmat{\eta_2 \cr \phi_2} 
= \pmat{\zst\eta_1 \cr \zsh\phi_1}  .
\end{equation}
Noting that
$\zst\D \hst = \D\zsh h'\hsh = \D (1/\zeta') $, 
\qref{E.evp1} becomes
\begin{equation}\lbl{e.evp2}
\pmat{ 
\lam-\D \pfrac{1-u_1}{\zeta'} 
& -D\tanh D\cr  
\tig \pfrac{1-v_1}{\zeta'}  &
\lam-\pfrac{1-u_1}{\zeta'}\D} 
\pmat{\eta_2 \cr \phi_2} 
=0,
\end{equation}
where (with formal leading order behavior indicated)
\begin{equation}
 u_1 = \zsh u \sim \uu, 
\qquad 
v_1 = \zsh\tig\inv(1-u)v' =\frac{(1-\zsh u)\D\zsh v}{\tig(1+\uu)}
\sim \zsh v'.
\end{equation}

\noindent
{\bf Approximate diagonalization.}
In order to reduce the eigenvalue problem to the `right' scalar equation,
it is helpful to balance off-diagonal terms (up to a commutator) 
and diagonalize the leading part of the operator. 
Let us define $p$, $q$, $u_p$, $u_q$,
and for later reference also $\rho$ and $u_\rho$,
so that
\begin{equation}\lbl{d.pq}
{\frac{1-u_1}{\zeta'}} = p = 1+u_p,
\qquad 
 \sqrt{\frac{1-v_1}{\zeta'}} = q = 1+u_q,
\qquad \sqrt q=\rho=1+u_\rho. 
\end{equation}
Asymptotically we expect 
\begin{equation}\lbl{e.upuq}
u_p \sim -2\uu, \qquad  u_q \sim -\sfrac12 \uu
\qquad u_\rho\sim-\sfrac14 \uu.
\end{equation}
To make precise estimates, we write 
\[
u_p(x)=\eph^2 \ti u_p(\eph x),\qquad
u_q(x)=\eph^2 \ti u_q(\eph x),  \qquad 
u_\rho(x)=\eph^2 \ti u_\rho(\eph x), 
\]
and apply the scaled $H^2$ bounds from Theorem~\ref{t.sw}
and \qref{i.taus} to the expressions in \qref{e.uvomeg}, 
using standard calculus inequalities. Straightforward
computations yield the following.
\begin{lemma}\lbl{l.coefs} 
For $\eph>0$ sufficiently small, 
the $H^2$ norms of $\ti u_p$, $\ti u_q$ and $\ti u_\rho$
are bounded by a constant $K$ independent of $\eph$, and 
the functions $u_p$, $u_q$, $u_\rho$ satisfy the pointwise bounds 
\begin{equation}\lbl{i.upuq}
|u_p| + |u_q|+|u_\rho| \le K\eph^2, 
\qquad 
|u_p'|+ |u_q'|+|u_\rho'|\le K\eph^{3}.
\end{equation}
Furthermore, as $\eph\to0$ we have
\begin{equation}\lbl{l.tiup}
\|\ti u_p + 2\sw\|_{H^1} \to 0 , \qquad
\|\ti u_q + \half\sw\|_{H^1} \to 0 .
\end{equation}
\end{lemma}

Introduce the operator (Fourier multiplier)
\begin{equation}\lbl{d.BBpm}
\SC= \sqrt{-\gamma D\tanh D}. 
\end{equation}
In order to balance orders of differentiation 
in the system, we change variables via 
\begin{equation}\lbl{d.eu3}
\pmat{\eta_3\cr \phi_3} = 
\pmat{\tig q \eta_2\cr \SC \phi_2} .
\end{equation}
The system \qref{e.evp2} then takes the (partially symmetrized) form
\begin{equation}\lbl{e.ev3}
\pmat{
\lam-\D p + R_1 &  q\SC \cr
\SC q & \lam - \D p + R_2
}
\pmat{\eta_3\cr \phi_3} =  0, 
\end{equation}
where $R_1$ and $R_2$ (which will both turn out to be negligible)
are given by
\begin{eqnarray}
&& R_1 = \D p - q\D pq\inv = (\D q-q\D)pq\inv =q' pq\inv ,
\lbl{d.r1}\\
&& R_2 = \D p - \SC p\D\SC\inv = 
p'+[p,\SC]\SC\inv\D. 
\lbl{d.r2}
\end{eqnarray}
Finally, we approximately diagonalize by changing variables via
\begin{equation}\lbl{d.4}
\pmat{\eta_4\cr \phi_4}  =
\pmat{1 &-1\cr 1 & 1}
\pmat{\eta_3\cr \phi_3}  , 
\qquad
\frac12
\pmat{1 &1\cr -1 & 1}
\pmat{\eta_4\cr \phi_4}  =
\pmat{\eta_3\cr \phi_3}  . 
\end{equation}
Then the system \qref{e.ev3} takes the form
\begin{equation}\lbl{e.ev4b}
\left(
\pmat{\lam-\D p &0\cr 0&\lam - \D p} +
\frac12\pmat{
-\SC q-q\SC +R_1+R_2 &
-\SC q+q\SC +R_1-R_2 \cr
\SC q-q\SC +R_1-R_2 &
\SC q+q\SC +R_1+R_2 
}
\right)
\pmat{\eta_4\cr \phi_4}  = 0.
\end{equation}

We make a few observations regarding the form of this system:
First, the operators $R_2$ and $[\SC,q]=\SC q-q\SC$
involve commutators and will turn out to be bounded, 
while $R_1$ is just a multiplier. 
So the off-diagonal terms are bounded operators, involving no derivatives.
Second, as is needed for energy estimates, we will invoke the symmetrization
identity 
\begin{equation}\lbl{e.sym1}
\D p=\sqrt p\,\D \sqrt p+ \half p' .
\end{equation}
The operator 
$\frac12(\SC q+q\SC)$ 
can be explicitly symmetrized (for energy estimates)
up to a (double) commutator in terms of $\rho=\sqrt q$:
\begin{equation}\lbl{e.sym2}
\sfrac12(\SC q+q\SC) 
= \sqrt q\,\SC\sqrt q + \sfrac12[[\SC,\sqrt q],\sqrt q].
\end{equation}
Finally, note that the weight-transformed operator 
$\eax\SC\emax$ is a Fourier multiplier with symbol
\begin{equation}\lbl{e.DBBsym}
\SC(k+ia) = 
\sqrt{-\gamma \xi\tanh \xi}, \qquad \xi=k+ia.
\end{equation}
The principal square root is used here and
the real part is nonnegative. We define
\begin{equation}\lbl{d.Apm}
\AAp=\D+\SC, \qquad \AAm=\D-\SC.
\end{equation}
It is easy to see that these formulae define closed
operators in $L^2_a$ with domain $H^1_a$ 
and with 
spectrum given by the range of the weight-transformed Fourier
multipliers
\[
k\mapsto \Ac_{\pm}(\xi) = i\xi\pm \sqrt{-\gamma\xi\tanh\xi},
\qquad \xi=k+ia, \quad k\in\R.
\]

\noindent
{\bf Final form as system.}
Based on these observations, it will be 
convenient to write the eigenvalue problem as follows.
We use \qref{d.pq} to write the operator 
in the (1,1) and (2,2) slots of \qref{e.ev4b} 
as  $\lam-\Ac_{11}$ and $\lam-\Ac_{22}$
respectively, with
\begin{equation}\lbl{d.A11}
\Ac_{11} = 
\AAp+\UU + J_{11}, 
\qquad
\UU := \D u_p+\SC u_q,
\end{equation}
\begin{equation}\lbl{e.Apm}
\Ac_{22} = \AAsm + J_{22},
\qquad \BB_\pm := \sqrt p\,\D \sqrt p \pm \sqrt q\,\SC\sqrt q.
\end{equation}
The system \qref{e.ev4b} then takes the form
\begin{equation}\lbl{e.evfinal}
(\lam-\Ac)
\pmat{\eta_4\cr \phi_4}   
= 0,
\qquad
\Ac = \pmat{\Ac_{11} & J_{12}\\ J_{21} & \Ac_{22}},
\end{equation}
with the `junk terms' $J_{ij}$ given in terms of
$R_1=q'pq\inv$ and  $R_2 = p'+[p,\SC]\SC\inv\D$
by
\begin{equation} \lbl{d.junkmat}
\pmat{ J_{11} & J_{12}\\ J_{21} & J_{22}} = 
-\frac12 \pmat{
R_2+R_1+[\SC,q] &
R_1-R_2-[\SC,q] \\
R_1-R_2+[\SC,q] &
R_1+[p,\SC]\SC\inv\D+[[\SC,\rho],\rho] 
}.
\end{equation}
Another way we will sometimes use to write the (1,1) component of $\Ac$ is 
\begin{equation}\lbl{d.A11b}
\Ac_{11} = \AAsp+\ti J_{11},
\qquad
\ti J_{11}=-\half(R_1+[p,\SC]\SC\inv\D-[[\SC,\rho],\rho]).
\end{equation}

Our main results to be proved in this paper now amount to the following.
\begin{theorem}\lbl{t.Amain} 
(Asymptotic stability with weights)
Fix $\aaw\in(0,\half]$ and set $a=\aaw\eph$. 
For $\eph>0$ sufficiently small the following hold: 
\begin{itemize}
\item[(i)]  With domain $(H^1_a)^2$, $\Ac$ is the generator of a 
$C^0$ semigroup on $(L^2_a)^2$. 
\item[(ii)]
Whenever $\re\lam\ge-\frac16\aaw\eph^3$ with $\lam\ne0$, 
$\lam$ is in the resolvent set of $\Ac$. 
\item[(iii)] The value $\lam=0$
is a discrete eigenvalue of $\Ac$ with algebraic multiplicity 2. 
\item[(iv)] Restricted to the $\Ac$-invariant spectral complement
$\bar Y_a$ of the generalized kernel of $\Ac$,
the semigroup $e^{At}$ is asymptotically stable, satisfying
\[
\|e^{\Ac t}z\|_a \le Ke^{-\beta t}\|z\|_a
\] 
for all $t\ge0$ and $z\in \bar Y_a$, with some constants $K>0$
and $\beta>\frac16\aaw\eph^3$ depending on $\eph$ and $\aaw$.
\end{itemize}
\end{theorem}
\begin{theorem}\lbl{t.nowt2} (Spectral stability without weights)
For $\eph>0$ sufficiently small, 
with domain $(H^1)^2$ in the space $(L^2)^2$, 
the spectrum of the operator $\Ac$ is precisely the imaginary axis.
\end{theorem}

{\bf Equivalences.} The statements in these correspond directly to those in Theorem~\ref{t.main}
due to the following facts. First, the map 
$(\eta_2,\phi_2)\mapsto(\eta_4,\phi_4)$
is clearly an isomorphism from $Z_a=L^2_a\times H^{1/2}_a$ to 
$(L^2_a)^2$ when $a>0$. (Note the symbol $\SC(k+ia)$ does not vanish
at $k=0$ in this case). Second, the map $\eta_1\mapsto\eta_2$ from \qref{d.eU2} 
is clearly an isomorphism on $L^2_a$ (and on $H^1_a$). 
For example,
\[
\|\eta_1\|_a^2 = 
\intR e^{2as}|\eta_1\circ\zeta(s)\zeta'(s)|^2 e^{2a(\zeta(s)-s)}
\,\frac{ds}{\zeta'(s)}
=(1+O(\eph^2))\|\eta_2\|_a^2,
\]
since pointwise $|\zeta'-1|+\aaw\eph|\zeta(s)-s|=O(\eph^2)$ uniformly,
due to Theorem~\ref{t.sw}. Next, the composition map
$\zsh$ is an isomorphism on $H^s$ for $s=0$, 1 and 2, 
hence also for $s=\frac12$ and $\frac32$ by interpolation
(see \cite{BeLf}, particularly Theorems 3.1.2 and 6.4.4).
Therefore the map
\[
\phi_1\mapsto\phi_2=\emax e^{-a(\zeta(x)-x)}\zsh \eax\phi_1
\]
is an isomorphism on $H^s_a$.
(Note that the multipliers $e^{\pm a(\zeta(x)-x)}=I+ O(\eph^2)$ 
on $H^s$ for $s=0$, 1 and 2, as is easy to check 
using Theorem~\ref{t.sw}.) 

Finally, we claim that the transformation steps 
\qref{d.eU2}, \qref{d.eu3} and \qref{d.4} 
map the space of pairs $(\eta_1,\phi_1)$ satisfying condition \qref{e.espc} 
of Theorem~\ref{t.nowt}
isomorphically to the space of pairs $(\eta_4,\phi_4)\in 
L^2\times L^2$.
The first step to show this is to see that $(\eta_4,\phi_4)\in L^2\times L^2$
is equivalent to finiteness of the linearized energy: 
\begin{equation}\lbl{e.espc2}
\intR 
\phi_1(-\Dx\Heta)\phi_1+ \tig\eta_1^2 \, dx<\infty.
\end{equation}
The key point here is that since $\phi_2=\zsh\phi_1$ and $\Heta=\zshi\Hil\zsh$,
due to \qref{E.zrel} the change of variables $x=\zeta(\xf)$,
$dx=\zeta'(\xf)d\xf$ yields
\[
\intR \phi_1(-\Dx\Heta)\phi_1\,dx = 
\intR (\zsh\phi_1)(-\zsh \Dx\Heta \phi_1)\zeta'd\xf = 
\intR \phi_2 (D\tanh D)\phi_2\,d\xf = \|\phi_3\|_{L^2}^2.
\]
The second step is to demonstrate an equivalence of norms. We claim

\begin{lemma}
For some positive constants $c_-$ and $c_+$ 
independent of $\phi_1$,
\begin{equation}
c_- \intR \phi_1(D\tanh D)\phi_1 \, dx \le 
\intR \phi_1(-\Dx\Heta)\phi_1 \, dx \le
c_+ \intR \phi_1(D\tanh D)\phi_1 \, dx. 
\end{equation}
\end{lemma}
For the proof, it is enough to consider $\phi_1$ smooth and rapidly decaying on $\R$.
Let $\phi$ and $\phif$ be functions that are harmonic in
the fluid domain $\dometa$ and the flat strip $\domflat$ respectively,
and satisfy the boundary conditions
\begin{equation}\lbl{bc.phis}
\phi_1(x) = \phi(x,\eta(x))= \phif(x,0), 
\qquad 0=\phi_y(x,-1)= \phif_y(x,-1).
\end{equation}
Then 
\begin{equation}\lbl{e.phifl}
\intR \phi_1(-\Dx\Heta)\phi_1 \, dx 
= \intR\int_{-1}^{\eta(x)} |\nabla \phi|^2 \,dy\,dx,
\end{equation}
\begin{equation}\lbl{e.phif}
\intR \phi_1(D\tanh D)\phi_1 \, dx = \intR\int_{-1}^0 |\nabla \phif|^2 \,dy\,dx.
\end{equation}
Moreover, the function $\phi$ (resp. $\phif$) {\em minimizes} the double integral 
in \qref{e.phifl} (resp. \qref{e.phif}) among functions satisfying the 
same Dirichlet boundary conditions.
Let $X:\domflat\to\dometa$ be a smooth (but non-conformal) change of variables 
of the form $X(x,y) = (x,\tilde y(x,y))$, such that $X(x,0)=(x,\eta(x))$.
The function $\tilde \phi = \phi\circ X$ is smooth on $\domflat$ and satisfies
$\tilde\phi(x,0)=\phi(x,\eta(x))=\phi_1(x)$.
Using the minimizing property of $\phif$, then changing variables
and using that the gradient and (inverse) Jacobian of $X$ are uniformly bounded, 
we find
\[
\intR\int_{-1}^0 |\nabla\phif|^2\,dy\,dx \le
\intR\int_{-1}^0 |\nabla\tilde\phi|^2\,dy\,dx \le
\frac1{c_-}
 \intR\int_{-1}^{\eta(x)} |\nabla \phi|^2 \,dy\,dx .
\]
This establishes the first inequality in the Lemma.
The other one is similar.

\section{Estimates on commutators and junk}\lbl{S.comm}

In order to bound the junk terms, we need to bound commutators of
$\SC$ with the multipliers $p$, $q$ and $\rho$, or equivalently with
$u_p$, $u_q$, and $u_\rho$, since $[\SC,1]=0$. 
The functions $u_p$, $u_q$, $u_\rho$ all have the 
scaled form $\eph^2 G(\eph x)$, where $G$ depends on $\eph$ but
remains bounded in $H^2$. 
The following result provides a general estimate for the commutator of
a Fourier multiplier and a multiplier with this scaled form.
Write $\ip{ k}^s=(1+k^2)^{s/2}$, $\ip{D}^s=(1+D^2)^{s/2}$.

\begin{prop}\lbl{P.comm}
Let $\PP$, $\QQ$ and $\RR$ be Fourier multipliers 
with symbols $P$, $Q$ and $R$ respectively, and  let $s\ge0$.  
Let $g(x)=\eph^2 G(\eph x)$ where $G\colon\R\to\R$ is smooth and exponentially
decaying, and let $f\colon\R\to\R$ be smooth with compact support.
Then
\[
\| \PP[\QQ,g]\RR f\|_{L^2} \le C_* C_G \|f\|_{L^2},
\]
where 
\[
C_* = \sup_{k,\hat k\in\R} \eph^2 \frac{
P(\eph k) |Q(\eph k)-Q(\eph\hat k)| R(\eph \hat k)
}{\ip{k-\hat k}^s},
\qquad 
C_G = \intR \ip{k}^s|\hat G(k)|\,\frac{dk}{2\pi}.
\]
\end{prop}

\noindent{\it Proof.} Using the Fourier transform and Young's
inequality, since $\hat g(k)=\eph\hat G(k/\eph)$, we have
\begin{eqnarray*}
\| \PP[\QQ,g]\RR f\|_{L^2}^2 &=& 
\intR \left| \intR P(k) (Q(k)-Q(\hat k))\eph 
\hat G\left(
\frac{k-\hat k}{\eph} 
\right)
R(\hat k) 
\hat f(\hat k)\,\frac{d\hat k}{2\pi}\right|^2 \frac{dk}{2\pi}
\\
&\le& C_*^2 \intR 
\left(\intR 
\ipbig{ \frac{k-\hat k}{\eph} }^s
\left|\hat G\left(
\frac{k-\hat k}{\eph} 
\right)\right| |\hat f(\hat k)|\,
\frac{d\hat k}{2\pi\eph}\right)^2\frac{dk}{2\pi}
\\
&\le& C_*^2 C_G^2 \|f\|_{L^2}^2.
\end{eqnarray*}
\begin{cor} \lbl{c.comm1}
Suppose $0\le a<\pi/4$.
With $g=u_p$, $u_q$ or $u_\rho$, there exists $K>0$ such that for all
small enough $\eph>0$, we have the $L^2_a$ operator norm estimates
\begin{equation}\lbl{i.DBg}
\|[\SC,g]\SC\inv\D\|_a \le K\eph^3,
\end{equation}
\begin{equation}\lbl{i.jjjj}
\|J_{11}\|_a+ \|J_{12}\|_a+ \|J_{21}\|_a+ \|J_{22}\|_a \le K\eph^3.
\end{equation}
\end{cor}

\noindent{\it Proof.}
Observe that for each of the indicated choices for $g$, we have 
that $G$ is uniformly bounded in $H^2$ as a consequence of Lemma~\ref{l.coefs}.
So, using $s=2-\frac23$ we have
\begin{equation}\lbl{i.CG}
C_G\le \left(\intR \ip{k}^{-4/3}dk\right)^{1/2} 
\left(\intR \ip{k}^4| \hat G(k)|^2 dk\right)^{1/2} \le K
\end{equation}
independent of $\eph$.
And, the operator $\ip{D}^{-1/2}\SC\inv\D$ is uniformly bounded on $L^2_a$ since
its weight-transformed symbol is 
$\ip{\xi}^{-1/2}i\xi/ \sqrt{-\gamma\xi\tanh\xi}$, which is uniformly bounded.
Hence it suffices to show that with the choices $P(k)=1$,
$Q(k)=\sqrt{-\xi\tanh\xi}$, $R(k)=\ip{\xi}^{1/2}$ ($\xi=k+i\aw$), 
we have
\begin{equation}\lbl{i.CR}
C_*\le K\eph^3. 
\end{equation}
To prove this estimate the idea is to show that 
with $\xi=k+i\aw$, $\hat\xi=\hat k+i\aw$,
\begin{equation}\lbl{i.Rlip}
|Q(k)-Q(\hat k)| =
\left|\frac{\xi\tanh\xi-\hat\xi\tanh\hat\xi}
{Q(k)+Q(\hat k)}
\right|
\le \frac{K|k-\hat k|}{\max(1,|\hat \xi|^{1/2})}
\le \frac{K|k-\hat k|}{\ip{\hat \xi}^{1/2}},
\end{equation}
and conclude through scaling by $\eph$. 

To prove \qref{i.Rlip}, first note that 
$|Q(k)-Q(\hat k)|\le K|k-\hat k|$, 
since $Q'(k)$ is uniformly bounded, as is easy to show.
Suppose now that $\hat k>1$, without loss. If $k<0$
then
\[
|Q(k)-Q(\hat k)|\le
|Q(k)-Q(0)|+|Q(0)-Q(\hat k)|\le K(|k|+|\hat k|)=K|k-\hat k|.
\]
If $k>0$, then one computes explicitly that
\[
-\xi\tanh\xi = -(k+ia) \frac{\sinh 2k+i\sin 2a}{\cosh 2k+\cos 2a} ,
\]
and finds that $Q(k)$ lies in the fourth quadrant of the 
complex plane together with $Q(\hat k)$, and so 
\[
|Q(k)+Q(\hat k)|\ge |Q(\hat k)|\ge \ip{\hat \xi}^{1/2}/K,
\]
since $\tanh\hat \xi$ is bounded away from zero.
As the map $k\mapsto -\xi\tanh\xi$ is uniformly Lipschitz, the
estimate \qref{i.Rlip} follows. 
This proves the commutator estimate \qref{i.DBg}.
Using this together with \qref{i.upuq} and the fact that $\SC\D\inv$ 
is bounded, the remaining estimates in \qref{i.jjjj} follow directly.

\section{Symbol expansions and estimates}

Here we develop basic approximations and key estimates
that concern the Fourier multiplier $\Ac_+$.

\subsection{Low-frequency expansion and KdV scaling}
The Taylor expansion of $\tanh$ at zero,
\[
\tanh\xi = \xi-\sfrac13\xi^3+O(\xi^5),
\]
and the fact that for $\xi$ with positive imaginary part 
one has $\sqrt{-\xi^2}=-i\xi$, yields
\begin{equation}\lbl{e.expand1}
\sqrt{-\xi\tanh\xi} = -i\xi+\sfrac16i\xi^3+O(\xi^5).
\end{equation}
We note that if $\kxi\in\R$ and $|\kxi|$ is sufficiently small,
\begin{equation}\lbl{i.tanh2}
\sqrt{\frac{\tanh \kxi}{\kxi}} \le 1-\sfrac19\kxi^2.
\end{equation}
In the KdV long-wave scaling, one
replaces $\xi$ in \qref{e.expand1} by $\eph\xi$, and $\lam$ by $\eph^3\lah$.
We find (recall $\gamma=1-\eph^2$)
\begin{equation}\lbl{e.expand2}
\Ac_+(\eph\xi) = 
i\eph\xi+\sqrt{-\gamma\eph\xi\tanh\eph\xi} = 
 \eph^3( \sfrac12 i\xi\gamma_1 + \sfrac16 i\xi^3\gamma+\xi^3O(\eph^2\xi^2)),
\end{equation}
where 
\[
\gamma_1=2\eph^{-2}(1-\sqrt{1-\eph^2})=1+O(\eph^2).
\]
The KdV-scaled weight-transformed symbol of $-\lam+\AAp$ with $\xi=k+i\aaw$,
$\lam=\eph^3\lah=\eph^3(\lah_r+i\lah_i)$ is
\begin{align}
&\eph^{-3}(-\eph^3\lah+i\eph\xi+\sqrt{-\gamma\eph\xi\tanh\eph\xi}) = 
-\lah+\sfrac12 i\xi\gamma_1+\sfrac16i\xi^3\gamma
+\xi^3O(\eph^2\xi^2)
\nonumber \\
&\qquad = 
(-\lah_r -\sfrac12\aaw\gamma_1+ \sfrac16(\aaw^3-3\aaw k^2)\gamma)
\ + \  i( -\lah_i + \sfrac12 k\gamma_1 + \sfrac16(k^3-3k\aaw^2)\gamma )
+\xi^3O(\eph^2\xi^2) . 
\lbl{e.expand3}
\end{align}
This corresponds to the purely formal KdV approximation
(writing $\AAp=\AAp(D)$)
\[
-\lam+\AAp(\eph D) \sim \eph^3(-\lah+\half\D-\sfrac16\D^3).
\]

\subsection{High-frequency estimates}
\noindent
\begin{lemma}\lbl{L.sqrt} If $z\in\C$, $\aw>0$, then 
$\re\sqrt{z}\le \aw$ if and only if $\frac12(|z|+\re z) \le \aw^2$.
\end{lemma}

\smallskip\noindent
{\it Proof.} Write $\sqrt{z}=u+iv$
where $u\ge0$. The result follows from 
\[
|z|=|u+iv|^2=u^2+v^2, \quad \re z=u^2-v^2.
\]

\begin{lemma}\lbl{L.dispest}
Suppose $ \xi=\kxi+i\aw $ with $\kxi\in\R$ and $0<\aw<\pi/4$.
Then 
\begin{equation}\label{rew1a}
 0< \re \sqrt{-\xi\tanh\xi}\le 
\frac{a}{\sqrt{\cos2a}} 
\sqrt{\frac{\tanh \kxi}{\kxi}}.
\end{equation}
\end{lemma}

\smallskip
\noindent{\it Proof.} By the previous lemma, 
if $\beta>0$ and $w=\sqrt{-\xi\tanh\xi}$,
then $\re w\le \beta$ if and only if 
\begin{equation}\label{re.1}
|\xi||\tanh \xi| - \re(\xi\tanh\xi) \le 2\beta^2.
\end{equation}
We may write
\[
\tanh \xi = 
\frac {e^{\xi}-e^{-\xi}} {e^{\xi}+e^{-\xi}} = 
\frac{\sinh 2\kxi+i\sin 2\aw}
{\cosh 2\kxi+\cos 2\aw}
= \frac{u+iv}{D_0}
\]
with $u=\sinh2\kxi$, $v=\sin2\aw$, $D_0=
{\cosh 2\kxi+\cos 2\aw}.$
Then \qref{re.1} is equivalent to 
\[
(\kxi^2+\aw^2)^{1/2}(u^2+v^2)^{1/2}-(\kxi u-\aw v) \le
2\beta^2 D_0,
\]
or (taking $\kxi>0$ without loss)
\begin{equation}\label{re.2}
\left(1+\frac{\aw^2}{\kxi^2}\right)^{1/2}
\left(1+\frac{v^2}{u^2}\right)^{1/2} - 1+\frac{\aw v}{\kxi u}
\le 2\beta^2 \frac{D_0}{\kxi u}.
\end{equation}
Note that 
\[
\frac vu = \frac{\sin 2\aw}{\sinh 2\kxi} \le
\frac{\aw}{\kxi}. 
\]
Using this bound on the left-hand side of \qref{re.2}, we find that
\qref{re.2} is implied by the bound
\begin{equation}\lbl{re.2b}
\frac{\aw^2}{\kxi^2} 
\le \beta^2 \frac{D_0}{\kxi u}.
\end{equation}
Since $D_0\ge (\cosh2\kxi+1)\cos 2a = 2\cosh^2\kxi\cos2a$
and $u=2\cosh\kxi\sinh\kxi$, 
\qref{re.2b} is implied by 
\begin{equation}
\frac{a^2}{\cos 2a}\frac{\tanh\kxi}{\kxi}= \beta^2.
\end{equation}
This yields \qref{rew1a} as claimed.

\begin{cor}\label{c.DBsym}
Take $\gamma=1-\eph^2$, $\aw=\eph\aaw$ with
$0<\aaw\le\frac12$.  For $\eph>0$ sufficiently small,
we have that for all real $\kxi$,  with $\xi=\kxi+i\aw$,
\[
\re\sqrt{-\gamma\xi\tanh\xi} \le \eph\aaw
(1-\sfrac14\eph^2)
\sqrt{\frac{\tanh \kxi}{\kxi}} , 
\]
\[
\re\Ac_+(\xi)=
\re( i\xi+\sqrt{-\tig\xi\tanh\xi}) \le  \eph\aaw\left(-1+
(1-\sfrac14\eph^2)
\sqrt{\frac{\tanh \kxi}{\kxi}}\right) \le -\frac14\eph^3\aaw.  
\]
Moreover, uniformly for $\re\lam\ge-\frac16\eph^3\aaw$, the 
$L^2_a$ operator norm of the resolvent of $\Ac_+$
satisfies
\begin{equation}\lbl{i.AAp}
\|(\lam-\Ac_+)\inv\|_a\le \frac{12}{\aaw\eph^3}.
\end{equation}
\end{cor}
The first inequality follows since 
\[
\frac{\gamma}{\cos 2a }\le\frac{1-\eph^2}{1-2\eph^2\aaw^2}
\le 1-\sfrac12\eph^2 \le (1-\sfrac14\eph^2)^2,
\]
and the resolvent bound follows since 
$|\lam -\Ac_+(\xi)|\ge 
\re(\lam-\Ac_+(\xi))\ge \frac1{12}\eph^3\aaw$
for all $k$.

For later reference, we note that for $\xi=\kxi+i\aw$ with 
$|a|<\pi/8$, $\kxi\in\R$,
\begin{equation}\lbl{i.tanh}
|\tanh\xi| \le 
\frac{ |\sinh 2\kxi|+|\sin 2a|}{\cosh 2\kxi+\cos 2a} \le 1.
\end{equation}

%
%
\section{Semigroup generation and scalar reduction by elimination}

To start our analysis of the linearized dynamics governed by $\Ac$, 
we use energy estimates to establish resolvent bounds for the symmetrized
operators $\AAspm$ that dominate the diagonal of $\Ac$.
By consequence, we show in this section that $\Ac$ generates
a $C^0$ semigroup in $(L^2_a)^2$, with $a=\eph\aaw$ for $\aaw\in[0,\half]$. 
Also we will show that if $\aaw\in(0,\half]$, 
 $\lam-\Ac_{22}$ in \qref{e.evfinal} is uniformly invertible 
on $L^2_a$ for all $\lam$ satisfying $\re\lam\ge-\half\eph\aaw$.
This allows us to eliminate $\phi_4$ in the eigenvalue problem 
\qref{e.evfinal} and reduce to a scalar,
nonlinear eigenvalue equation for $\eta_4$ in the form
\begin{equation}\lbl{e.eve1}
(\lam-\Ac_{11}- J_{12}(\lam-\Ac_{22})\inv J_{21} )\eta_4 = 0.
\end{equation}


\begin{lemma}\lbl{l.AAspm}
For some constant $K$ independent of $\eps$, $\aaw$ and $\lam$,
if $\aaw\in[0,\half]$, $\eps>0$ is sufficiently small, 
and $\re\lam>-\eph\aaw(1-K\eph^2)$,
then $\lam$ is in the resolvent set of $\AAsm$,
with
\begin{equation}\lbl{i.AAmi}
\|(\lam-\AAsm)\inv\|_a \le \frac{1}{\re\lam + \aaw\eph(1-K\eph^2)}, 
\end{equation}
and if $\re\lam> K\aaw\eph^3$ then $\lam$ is in the resolvent set of $\AAsp$,
with 
\begin{equation}\lbl{i.AApi}
\|(\lam-\AAsp)\inv\|_a \le \frac{1}{\re\lam - K\aaw\eph^3}. 
\end{equation}
\end{lemma}

{\it Proof.} The main step in the proof is to perform energy estimates 
for each term. Let $z$ be smooth with compact support, 
and let $z_a=\eax z$. Write $\D_a=\eax\D\emax=\D-a$, and 
recall $a=\alpha\eph$. 
We compute (using the fact that $p\ge 1-K\eph^2$ in the last step)
\begin{eqnarray}
\re\ip{-\sqrt p\D\sqrt p \, z,z}_a &=& 
-\re\intR (\sqrt p\D\sqrt p \, z) \bar z\, e^{2ax}dx 
= -\re\intR (\D_a\sqrt p\, z_a)\overline{\sqrt pz_a}\,dx
\nonumber\\
&=& a\intR p|z_a|^2\,dx \ge \alpha\eph(1-K\eph^2)\|z\|_a^2.
\label{i.pz}
\end{eqnarray}
%
Due to \qref{e.DBBsym} above, with $\xi=k+i\aw$,
\begin{equation}
\re\ip{\sqrt q\SC\sqrt q\, z,z}_a =
\re\intR (\sqrt q\SC\sqrt q z_a)\bar z_a e^{2ax}\,dx
= \re\intR 
\sqrt{-\gamma\xi\tanh\xi}|\FF \sqrt q\, z_a|^2 \,\frac{dk}{2\pi}.
\end{equation}
By Corollary~\ref{c.DBsym} we find $0\le 
\re\ip{\sqrt q\SC\sqrt q\, z,z}_a \le \eph\aaw\|qz\|_a^2\le \eph \aaw(1+K\eph^2)\|z\|_a^2.$
Hence it follows that for all smooth $z$ with
compact support
\begin{eqnarray}
\frac{\re\ip{(\lam-\AAsm)z,z}_a}{\|z\|_a^2} &\ge& \re\lam+ \alpha\eph(1-K\eph^2), 
\lbl{i.Asm}
\\
\frac{\re\ip{(\lam-\AAsp)z,z}_a}{\|z\|_a^2} &\ge& \re\lam- K\alpha\eph^3 .
\lbl{i.Asp}
\end{eqnarray}
When the right-hand side is positive, this proves $\lam-\AApm$ is 
uniformly invertible {\it on its range}, satisfying the respective estimates in
\qref{i.AAmi} and \qref{i.AApi}.

To prove that $\lambda$ is in the resolvent set of $\AAsm$, 
what remains to prove is that the range of $\lam-\AAspm$ is all of $L^2_a$. 
To accomplish this, 
we use a perturbation estimate to establish that a fixed value
$\lam=1$ is in the resolvent set for small enough $\eps$,
then invoke an analytic continuation property of resolvents. 
For $\lam=1>0$ fixed, if $\eps$ is small then we will show
$1-\AAspm$ is a small relative perturbation of the Fourier multiplier
$1-\AApm$ from \qref{d.Apm},
with 
\begin{equation}\lbl{i.ptbA}
\|(1-\AApm)\inv(\AAspm-\AApm)\|_a\le K\eph^2<1.
\end{equation}
By perturbation it follows $1$ is in 
the resolvent set of operator $\AAspm$
and the range of $\lam-\AAspm$ is all
of $L^2_a$  for $\lam=1$. 
Using the Neumann series for the resolvent, 
we see that the resolvent of any closed operator can be 
analytically continued to any set
where the resolvent has a uniform apriori bound
(see Theorem~III.6.7 of \cite{Kato}). 
By consequence, the resolvent set of
$\AAsm$ (resp. $\AAsp$) contains the entire right half-plane where
the right-hand side of \qref{i.Asm} (resp. \qref{i.Asp}) is positive.

We proceed to prove \qref{i.ptbA}. We compute that 
\[
\AAspm-\AApm = \D u_p  -\half p' \pm \SC u_q \pm[\rho,\SC]\rho.
\] 
Since $\SC\D\inv$ is bounded, Corollary~\ref{c.comm1} and \qref{i.upuq} imply
\begin{equation}\lbl{i.cbdrho}
|u_p|+|u_q|+ |p'|+ \|[\SC,\rho]\rho\|_a 
\le K\eph^2.
\end{equation}
We claim that for some constant $K$ independent of $\eps$ and $\aaw$,
\begin{equation}\lbl{i.Amrbd}
\|(1-\AApm)\inv \D^j\|_a \le K
\quad\mbox{for $j=0$ and 1}. 
\end{equation}
The symbol of the weight-transformed operator 
$-\emax(1-\AApm)\inv \D^j\eax$ is 
\begin{equation}
m_j(\xi) = \frac{(i\xi)^j}{1-i\xi\mp\sqrt{-\gamma\xi\tanh\xi}},
\qquad \xi=k+i\aw,
\end{equation}
for $j=0$ or 1.
Since $-i\xi=-ik+a$, the real part of the denominator is always 
greater than 1, by Corollary~\ref{c.DBsym}.
Hence $|m_0(\xi)|\le1$ for all $k\in\R$, 
so \qref{i.Amrbd} holds for $j=0$.
For $j=1$, we have that $|m_1(\xi)|\le 9$ for $|\xi|\le9$,
while for $|\xi|>9$ the denominator is bounded below by
\[
|1-i\xi|-|\xi|^{1/2} \ge |\xi|(1-|\xi|^{-1/2})\ge \sfrac23|\xi|,
\]
because $|\tanh\xi|\le1$ by \qref{i.tanh}.
Thus $|m_1(\xi)|\le \frac32$ for $|\xi|\ge9$. 
Hence \qref{i.Amrbd} holds also for $j=1$.
The bound in \qref{i.ptbA} follows by combining
\qref{i.Amrbd} with \qref{i.cbdrho}.
This completes the proof of the Lemma.

\begin{prop}\lbl{p.semi} 
For $\eph>0$ sufficiently small, $\Ac$ is the generator of a
$C^0$ semigroup on $(L^2_a)^2$.
\end{prop}

{\it Proof.} By the Lemma just proved and the Hille-Yosida theorem, 
the operator
\[
\Ac_* = \pmat{ \AAsp & 0 \\ 0 & \AAsm}
\]
is the generator of a $C^0$ semigroup on $(L^2_a)^2$.
But $\Ac-\Ac_*$ is bounded, so the result follows from a 
standard perturbation theorem (see Theorem IX.2.1 of \cite{Kato}).

\begin{lemma}\lbl{p.Ac22} 
For some constant $K$ independent of $\eps$, $\aaw$ and $\lam$
and for $\aaw\in(0,\half]$, if $\eph>0$ is sufficiently small then
$\lam$ is in the resolvent set of $\Ac_{22}$
whenever $\re\lam\ge-\half\eph\aaw$, with
\begin{equation}\lbl{i.Ac22}
\|(\lam-\Ac_{22})\inv\|_a \le \frac{K}{\eph\aaw}.
\end{equation}
\end{lemma}

{\it Proof.} Due to the bound on $J_{22}$ from Corollary~\ref{c.comm1}, 
this result follows directly from the results in Lemma~\ref{l.AAspm}
concerning the resolvent of $\AAsm$.

\section{Resolvent bounds for $|\lambda|$ not too small}\lbl{S.large}

For the remainder of this paper we fix $\aaw$ satisfying
$0<\aaw\le\half$ and write $\aw=\eph\aaw$.
Here we demonstrate a bound on the resolvent
of the operator $\Ac$ from \qref{e.evfinal} that is uniform in $\lambda$,
for $\lambda$ in the right half-plane with $|\lambda|$ not too small.
This bound, in combination with 
the Gearhart-Pr\"uss spectral mapping theorem and our 
proof that the only eigenvalue of $\Ac$ in the right half-plane
is $\lam=0$ with algebraic multiplicity 2, will allow us to 
obtain linear asymptotic stability in $L^2_a$ for the semigroup $e^{\Ac t}$, 
conditional for perturbations containing no neutral-mode components.

\begin{prop}\lbl{p.resA}
Let $0<\nus<\nuh<1$.  For $\eph>0$ sufficiently small, 
all $\lambda$ satisfying 
\begin{equation}\lbl{c.large}
|\lam|\ge \eph^\nus \quad\mbox{and}\quad \re\lam\ge
-\sfrac1{12}{\aaw\eph^{1+2\nuh}}
\end{equation}
belong to the resolvent set of $\Ac$, with  
$\|(\lam-\Ac)\inv\|_a\le K/\eph^{1+2\nuh}$,
for some constant $K$ independent of $\eph$ and $\lam$.
\end{prop}


In the analysis we will make use of sharp Fourier cutoffs
(Fourier filters) defined as follows. We specify a wavenumber threshold
chosen to be $\kah=\eph^\nuh$, where we require $\nus<\nuh<1$. 
Define projection operators
$\pi_o$ (low-pass), $\pi_i$ (high-pass) on $L^2_a$  as follows.
First, on $L^2$, define low and high-pass filters by 
\begin{equation}\label{d.0lohi}
\pi_{0o}= \FF\inv\one_{[-\kah,\kah]} \FF,
\qquad \pi_{0i}= I-\pi_{0o}.
\end{equation} 
These operators are orthogonal projections on $L^2$. Now in $L^2_a$,
define orthogonal projections by  
\begin{equation}\label{d.lohi}
\pi_o = \emax\pi_{0o} \eax, \qquad \pi_i = I-\pi_o.
\end{equation}

\subsection{Resolvent bounds for $\Ac_{11}$ for $\lambda$ not small}
Crucial for our estimate of $(\lam-\Ac)\inv$ is to 
demonstrate uniform invertibility of the operator
\[
\lam-\D p-\SC q = \lam-\Ac_{11}-J_{11}
\]
which appears as the dominant part of the operator in
the $(1,1)$ slot of \qref{e.evfinal}. 
The aim here is to establish uniform invertibility
of the operator above
with respect to the weighted norm with weight $a=\alpha\eph$.
The estimate on the inverse will have the form
\begin{equation}\lbl{i.lampBq}
\|(\lam-\D p-\SC q)\inv \|_a \le \frac{K}{\eph^{1+2\nuh}} 
\end{equation}
and be valid for $\lam$ satisfying \qref{c.large},
provided $\eps>0$ is smaller than some fixed positive constant. 
(Here and below, $K$ is a generic constant independent of $\eps$
whose value may change from case to case.)
Since we know $\|J_{11}\|_a=O(\eph^3)$ by Corollary~\ref{c.comm1},
we infer that under
conditions of the same form on $\lam$ and $\eps$, $\lam-\Ac_{11}$ is
invertible with
\begin{equation}\lbl{i.A11}
\|(\lam-\Ac_{11})\inv\|_a\le \frac{K}{\eph^{1+2\nuh}}.
\end{equation}

To prove the bound \qref{i.lampBq} we study the equation 
\[
(\lam-\D p-\SC q)z = g
\]
decomposing this equation in terms of $z_o=\pi_oz$, $z_i=\pi_iz$,
$g_o=\pi_og$, and $g_i=\pi_ig$. Apply $\pi_o$ and note $\pi_oz_i=0$,
and the low-pass filter (nontrivially) commutes with derivatives 
and Fourier multipliers.  We get
\begin{equation}\lbl{e.zozi}
\pmat{\Ac_{oo} & \Ac_{oi}\cr
\Ac_{io} & \Ac_{ii}}
\pmat{z_o \cr z_i} =
\pmat{g_o \cr g_i}, 
\end{equation}
\[
\Ac_{oo} = 
\lam -\pi_o(\D p + \SC q)\pi_o,
\qquad 
\Ac_{oi} = -\pi_o(\D u_p+\SC u_q)\pi_i,
\] 
\[
\Ac_{io} = -\pi_i(\D u_p+\SC u_q)\pi_o,
\qquad
\Ac_{ii} = \lam -\pi_i(\D p + \SC q)\pi_i.
\]
Here recall $u_p= p-1$, $u_q=q-1$ satisfy the pointwise bounds
in \qref{i.upuq}.
The low-pass Fourier filter satisfies
$\|\pi_o\D\|_a\le |\kah+ia|\le 2\eph^\nuh$ and
$\|\pi_o\SC\|_a\le K_0\eph^\nuh$
since $|\xi\tanh\xi|^{1/2}=|\xi||\tanh\xi/\xi|^{1/2}
 \le K_0\eph^\nuh$ for $\xi=k+i\aw$ with $|k|\le\kah$,
with some constant $K_0$ independent of $\eph$. 
Also $\|p\pi_o\|_a+\|q\pi_o\|_a\le K_0$.
Now clearly, if 
$|\lambda|\ge \eph^\nus$ and $\eph$ is small enough 
so $\eph^\nus> 2K_1\eph^\nuh$ with $K_1=2K_0+K_0^2$
(we use $\nus<\nuh$ here),
then $\Ac_{oo}$ is invertible and
%
\begin{equation}\lbl{i.Aooinv}
\|\Ac_{oo}\inv\|_a \le \frac{1}{|\lam|-K_1\eph^\nuh} \le
\frac{2}{\eph^\nus}. 
\end{equation}
Since $\D \tilde u\pi_o=(\tilde u'+\tilde u\D)\pi_o$ 
for $\tilde u=u_p$ and $u_q$, we also find (since $\hat\nu\le1$)
\begin{equation}\lbl{i.Aoiio}
\|\Ac_{oi}\|_a\le K\eph^{2+\nuh}, 
\qquad
\|\Ac_{io}\|_a\le K\eph^{2+\nuh}. 
\end{equation}

In order to establish the estimate \qref{i.A11}, it suffices to show
that whenever $\lam$ satisfies \qref{c.large},
$\Ac_{ii}$ is invertible on $L^2_a$ with 
\begin{equation}\lbl{i.Aiiinv}
\|\Ac_{ii}\inv\|_a\le \frac{K}{\eph^{1+2\nuh}}
\end{equation}
This is because, after elimination of $z_0$, \qref{e.zozi} reduces to
\begin{equation}\lbl{e.zi}
(I - \Ac_{ii}\inv\Ac_{io}\Ac_{oo}\inv\Ac_{oi})z_i
= \Ac_{ii}\inv(g_i-\Ac_{io}\Ac_{oo}\inv g_o).
\end{equation}
Since
\[
\|\Ac_{ii}\inv\Ac_{io}\Ac_{oo}\inv\Ac_{oi}\|_a\le 
\frac{K}{\eph^{1+2\nuh}} (K\eph^{2+\nuh})^2 \frac{K}{\eph^\nus}
\le K\eph^{3-\nus},
\]
the desired estimate \qref{i.lampBq} 
then follows for sufficiently small $\eps$ 
and large $|\lam|$ from
\[
\|z_o+z_i\|_a \le \|z_o\|_a+\|z_i\|_a \le \frac{K}{\eph^{1+2\nuh}}
(\|g_i\|_a+\|g_o\|_a)
\le \frac{2K}{\eph^{1+2\nuh}} \|g\|_a.
\]

\subsection{Uniform invertibility of $\Ac_{ii}$ by energy estimates}
To prove the invertibility of $\Ac_{ii}$ with the estimate
\qref{i.Aiiinv},
since $z=z_o+z_i$, the main step is to prove the energy estimate
\begin{equation}\lbl{i.Aii1}
-\re\ip{\pi_i(\D p+\SC q)\pi_i z,z}_a \ge 
\sfrac1{10}{\eph^{1+2\nuh}\aaw}\|z_i\|_a^2
\end{equation}
for all smooth $z$ with compact support.
Given such $z$, let $z_a=\eax z$. Recall $\D_a=\eax\D\emax=\D-a$.
Then $\pi_{0i}z_a$ is in $H^m$ for all $m$ since
\[
\int_{|k|>\kah} (1+|k|^2)^m |\hat z_a|^2\,dk<\infty.
\]

1.  Recall $\D p=\sqrt p\D\sqrt p+\half p'$ and $|p'|\le K\eph^3$.
Similarly as for $\AAspm$ which we treated before, 
\begin{eqnarray}
-\re\ip{\pi_i\sqrt p\D\sqrt p \pi_i z,z}_a &=& 
-\re\intR (\pi_i\sqrt p\D\sqrt p \pi_i z) \bar z\, e^{2ax}dx 
\nonumber\\
&=& -\re\intR (\D_a\sqrt p \pi_{0i}z_a)\overline{\sqrt p\pi_{0i}z_a}\,dx
\nonumber\\
&=& a\intR p|\pi_{0i}z_a|^2\,dx \ge \eph\aaw(1-K\eph^2)\|\pi_i z\|_a^2.
\label{i.pz1}
\end{eqnarray}
The last inequality holds since $a=\eph\aaw$ and $p\ge1-K\eph^2$.

2. Now write $\rho=\sqrt q$ (as before) and compute
\begin{equation}
\pi_i\SC q \pi_i z
= \rho \pi_i\SC \rho z_i+ \CC_{\rho i}\rho z_i,
\end{equation}
where, with $u_\rho=\rho -1$ ($=O(\eph^2)$), we can write
$\CC_{\rho i} = \CC_\rho -\CC_{\rho o}$ with  
\begin{equation}
\CC_\rho = [\SC,\rho]=\SC u_\rho - u_\rho\SC ,
\qquad
\CC_{\rho o}= \pi_o\SC u_\rho-u_\rho\pi_o\SC, 
\end{equation}
Now, since $\|\pi_o\SC\|_a\le K\eph^\nuh$, 
evidently $\|\CC_{\rho o}\|_a \le K\eph^{2+\nuh}$, 
and Corollary~\ref{c.comm1} implies $\|\CC_\rho\|_a\le K\eph^3$.
Recall that the weight-transformed operator $\eax\SC\emax$ 
is a Fourier multiplier with symbol given by \qref{e.DBBsym}. 
For this 
we will use the high-frequency dispersion estimate in Corollary~\ref{c.DBsym}.
Note that for $|k|\ge \kah=\eph^\nuh$, if $\eph$ is small enough then
by \qref{i.tanh2},
\begin{equation}\lbl{d.hatdel}
1-\sfrac19\eph^{2\nuh}> 
\sqrt{\frac{\tanh \kah}{\kah}} \ge
\sqrt{\frac{\tanh k}{k}}.
\end{equation}
Using this with Corollary~\ref{c.DBsym}, we find
\begin{eqnarray}
-\re\ip{\rho\pi_i\SC \rho z_i,z_i}_a &=&
-\re\intR  (\rho \pi_i\SC \rho z_i)\bar z_i\, e^{2ax}dx
\nonumber \\
&=& -\re\int_{|k|>\kah} 
\sqrt{-\gamma\xi\tanh\xi} 
|\FF\rho z_{ia}|^2\,\frac{dk}{2\pi}
\nonumber \\
&\ge& -\eph \aaw\hat (1-\sfrac19\eph^{2\nuh})(1+K\eph^2) \|z_i\|_a^2.
\lbl{i.rhoDBi}
\end{eqnarray}

3. Combining \qref{i.pz1} with \qref{i.rhoDBi} and 
$|p'|+\|\CC_{\rho i}\rho\|_a\le K\eph^3$ yields \qref{i.Aii1}, since
for small $\eph$,
\[
\frac{-\re\ip{\pi_i(\D p+\SC q)\pi_i z,z}_a}{\|z_i\|_a^2} \ge 
\eph\aaw\left(1-K\eph^2 -(1-\sfrac19\eph^{2\nuh})(1+K\eph^2)
\right)-K\eph^3
\ge \sfrac1{10}\eph^{1+2\nuh}\aaw.
\]
Since $\Ac_{ii}=\lam\pi_o+\Ac_{ii}\pi_i$,
it follows that if $\lam$ satisfies \qref{c.large}, 
then $\Ac_{ii}$ has bounded inverse on its range
with bound given by \qref{i.Aiiinv}

4. To prove that the range of $\Ac_{ii}$ is all of $L^2_a$, we use the
same continuation approach as previously given for $\AAspm$. 
We can write
\[
\Ac_{ii} 
= \lam-(\D+\SC)\pi_i - \pi_i(\D u_p+\SC u_q)\pi_i
.  \]
For $\lambda=1$ fixed, the $\eps$-independent bound
\begin{equation}\lbl{i.relbd2}
\| (1-(\D+\SC)\pi_i)\inv \D\pi_i \|_a \le K
\end{equation}
follows by restricting the proof of \qref{i.Amrbd} to frequencies
$|k|\ge \kah$. Then we obtain
\[
\| (1-(\D+\SC)\pi_i)\inv \pi_i(\D u_p+\SC u_q) \|_a \le K\eph^2<1
\]
for small enough $\eps$, and  the invertibility of $\Ac_{ii}$ whenever
$\re\lam\ge0$ and $|\lam|\ge \eph/K$ now follows as before for
$\AAspm$, by continuation based on the energy estimate in step 3.

\subsection{Bound on the resolvent of $\Ac$}

Now we complete the proof of Proposition~\ref{p.resA}.
We solve the resolvent equation
\[
(\lam-\Ac)
\pmat{\eta_4 \cr \phi_4} =
\pmat{g_1 \cr g_2} 
\]
by simple elimination, writing 
\begin{eqnarray*}
\phi_4 &=& 
(\lam-\Ac_{22})\inv(g_2+J_{21}\eta_4),
\\
\eta_4 &=& W_*(\lam)\inv (\lam-\Ac_{11})\inv 
(g_1+ J_{12}(\lam-\Ac_{22})\inv g_2),
\\
W_*(\lam)&=& I-(\lam-\Ac_{11})\inv J_{12}(\lam-\Ac_{22})\inv J_{21},
\end{eqnarray*}
This is justified based on the estimates \qref{i.Ac22}, \qref{i.A11}, 
and the estimate
\begin{equation}\lbl{i.j1221}
\| J_{12}(\lam-\Ac_{22})\inv J_{21}\|_a \le K \eph^5
\end{equation}
that follows from Corollary~\ref{c.comm1} together with \qref{i.Ac22} 
for $\re\lam\ge -\half\eph\aaw$.
For the solution of the system, one obtains the estimates
\[
\|\eta_4\|_a\le \frac{K}{\eph^{1+2\nuh}}
\Big(\|g_1\|_a + K\eph^2\|g_2\|_a\Big), \qquad 
\|\phi_4\|_a\le \frac{K}{\eph}
\Big( \|g_2\|_a+ K\eph^{2-2\nuh}\|g_1\|_a \Big),
\]
whence the estimate $\|(\lam-\Ac)\inv\|_a\le K/\eph^{1+2\nuh}$
follows.

\section{KdV scaling and bundle limit}

It remains to study the eigenvalue problem when $|\lambda|$ is small.
satisfying $|\lam|\le\eph^\nus$.
At this point we have shown that the eigenvalue system \qref{e.evfinal} 
can be reduced to the nonlinear eigenvalue equation \qref{e.eve1} 
whenever $\re\lam\ge-\half\eph\aaw$. 
For $\re\lam\ge-\frac16\eph^3\aaw$,
we may further apply the Fourier multiplier
$(\lam-\AAp)\inv$ to \qref{e.eve1}, by Corollary~\ref{c.DBsym}. 
This reduces the eigenvalue problem to the 
nonlinear eigenvalue equation
\begin{equation}\lbl{d.W}
W(\lambda)\eta_4:= (I-(\lam-\AAp)\inv\UU - J_*)\eta_4 = 0,
\end{equation}
where 
\[
J_* = (\lam-\AAp)\inv(J_{11}+J_{12}(\lam-\Ac_{22})\inv J_{21}).
\]
The operator $J_*$ will be shown to be negligible.

As we shall see, the bundle $W(\lam)$ becomes singular at $\lam=0$
due to the fact that zero is an eigenvalue of the operator $\Ac$.
To determine the multiplicity of this eigenvalue and 
establish the invertibility of $W(\lam)$ for nonzero $\lam$, 
we  make use of the KdV long-wave scaling. 
We introduce scaled variables with tildes via 
\begin{equation}\lbl{d.kscal}
\ti x = \eph x, \quad \lam = \eph^3\ti \lam.
\end{equation}
Then $\D=\eph\ti\D$, and in purely formal terms we have the following
leading behavior (see \qref{e.expand2} and note $\sqrt{-D\tanh D}\sim -\D$):
\[
\lam-\AAp \sim \eph^3(\ti\lam-\half\ti\D+\sfrac16\ti\D^3),
\qquad \UU =\D u_p+\SC u_q \sim \eph^3 \ti\D (-2\sw)-\eph^3\ti\D(-\half\sw).
\]
Thus we expect $W(\lam)\sim W_0(\ti\lam)$ where (with tildes omitted on derivatives)
\begin{equation}\lbl{d.W0}
W_0(\ti\lam)= I+(\ti\lam-\half\D+\sfrac16 \D^3)\inv\D\left(\sfrac32\sw \right).
\end{equation}
This bundle $W_0(\ti\lam)$ is associated with the eigenvalue problem for 
the KdV equation scaled as
\[
\Dt f-\half\Dx f+ \sfrac32 f\Dx f+ \sfrac16\Dx^3 f =0,
\]
linearized about the soliton profile $f=\sw=\sech^2(\sqrt3x/2)$.

To be clear, what we are really doing when changing variables
is using a similarity transform in terms of the dilation operator $\Dil$
defined by
\begin{equation}\lbl{d.Se}
(\Dil f)(x)= f(x/\eph)/\sqrt{\eph}
\end{equation}
which maps $L^2_\aw$ isometrically onto $L^2_\aaw$ since $\aw=\aaw \eph$:
\[
\intR |f(x)|^2 e^{2ax}\,dx = 
\intR |\Dil f(y)|^2 e^{2\aaw y}\,dy , \qquad y=\eph x.
\]
(Note that similarity transform does not change operator norms,
but $\Dil\D\Dil\inv=\eph\D$.)

The formal discussion above involves uncontrolled approximations in terms
of derivatives.  But this motivates the following rigorous statement in
terms of convergence of bundles.
Based on this result, the scaled operator bundle will be studied using
the Gohberg-Sigal-Rouch\'e perturbation theorem \cite{GS71}.
\begin{prop}\lbl{t.wlim}
Define the scaled bundle $\tiW(\ti\lam) = \Dil W(\eph^3\ti\lambda)\Dil\inv$, 
and let
\begin{equation}\lbl{d.Omeh}
\Omeh: = \{\lah\in\C:
|\eph^3\lah|\le 1 , \ \re\lah\ge-\sfrac16\aaw\}.
\end{equation}
Then in operator norm on $L^2_\aaw$, we have
\begin{equation}\lbl{l.kdvop}
\sup_{\lah\in\Omeh}\|\tiW(\ti\lam)-W_0(\ti\lam)\|_\aaw \to 0 \quad\mbox{as $\eph\to0$}.
\end{equation}
\end{prop}
We prove this proposition by studying
pieces of the scaled bundle $\tiW(\ti\lam)$. 
Writing
$\Dil\SC\Dil\inv=\sqrt{-\tig\eph D\tanh\eph D}=\eph\ti\SC$, 
$\Dil\AAp\Dil\inv= \eph^3\ti\AAp$, and 
\[
u_p(x)=\eph^2 \ti u_p(\eph x),\qquad
u_q(x)=\eph^2 \ti u_q(\eph x),  \qquad 
\Dil J_*\Dil\inv =\ti J_*,
\]
the scaled bundle is written in the form
\begin{equation}\lbl{e.tiW}
\tiW(\ti\lam) = 
I-(\ti\lam-\ti\AAp)\inv(\D\ti u_p+\ti\SC\ti u_q) -\ti J_*.
\end{equation}
The proposition is implied by
following convergence results in operator norm on $L^2_\aaw$, 
to hold as $\eph\to0$, uniformly for $\ti\lam\in\Omeh$:
\begin{equation}\lbl{l.w1}
\| (\ti\lam-\ti\AAp)\inv\BB - (\ti\lam-\sfrac12\D+\sfrac16\D^3)\inv\D\|_\aaw \to 0
\quad\mbox{for both $\BB=\D$ and $\ti\SC$},
\end{equation}
\begin{equation}\lbl{l.w2}
\|\ti u_p + 2\sw\|_\aaw \to 0 , \qquad
\|\ti u_q + \half\sw\|_\aaw \to 0 ,
\end{equation}
\begin{equation}\lbl{l.w3}
\|J_*\|_a \to 0.
\end{equation}

In subsection~\ref{s.symlim}, 
we will prove the first limit \qref{l.w1} by studying 
the corresponding weight-transformed symbols. 
The third limit \qref{l.w3} is treated in subsection~\ref{s.jlim}.
The limits in \qref{l.w2} are a simple consequence of the fact that 
$\ti u_p+2\sw$ and $\ti u_q+\half\sw$ are pointwise multipliers, so
the $L^2_\aaw$ operator norm
is equal to the $L^\infty$ norm as a function,
and this is bounded by the $H^1$ norm,
which tends to zero by Lemma~\ref{l.coefs}.


\subsection{KdV limit for symbols}\lbl{s.symlim}

Here we establish the main limit \qref{l.w1} needed to prove the operator limit
in Theorem~\ref{t.wlim}.  Namely, we prove appropriate limits for
the scaled, weight-transformed symbol of the operator 
$\mesym(\lah,D)= (\lah-\ti\AAp)\inv \D$.
This symbol takes the form
\begin{equation}\lbl{d.mesym}
\mesym(\lah,\xi)= \frac{\eph^3 i\xi}
{-\eph^3\lah +i\eph\xi+\sqrt{-\gamma\eph\xi\tanh\eph\xi}},
\qquad \xi=k+i\aaw.
\end{equation}
The corresponding symbol for the limiting operator 
$\mosym(\lah,D)=(\lah-\half\D+\sfrac16\D^3)\inv\D$ is written
\begin{equation}\lbl{d.mosym}
\mosym(\lah,\xi) = \frac{ i\xi}{-\lah+\frac12i\xi+\frac16i\xi^3}.
\end{equation}
The symbol limits that we need to prove both limits in \qref{l.w1} are:
\begin{equation}\lbl{l.memo}
|\mesym(\lah,\xi)-\mosym(\lah,\xi)| \to 0
\qquad\mbox{as $\eps\to0$},
\end{equation}
\begin{equation}\lbl{l.memoB}
\left|\mesym(\lah,\xi)
\sqrt{\frac{\gamma\tanh\eph\xi}{\eph\xi}} -\mosym(\lah,\xi)\right| \to 0
\qquad\mbox{as $\eps\to0$}.
\end{equation}
These limits need to be established {\it uniformly for
$\xi\in\R+i\aaw$ and $\lah\in\Omeh$.}
(The second limit will follow easily once we establish the
first.)

1. First we provide simple, preliminary bounds on the limiting
symbol $\mosym$. As one sees from \qref{e.expand3}, 
for $\lah=\lah_r+i\lah_i$ and $\xi=k+i\aaw$,
the real part of the denominator of $\mosym$ is negative, 
with magnitude bounded below by
\begin{equation}\lbl{i.mo1}
\lah_r+\sfrac12\aaw(1-\sfrac13\aaw^2)+\sfrac12\aaw k^2
\ge \lah_r+\sfrac13\aaw+\sfrac12\aaw k^2
\ge \sfrac16\aaw(1+k^2)
\ge \sfrac16\aaw|\xi|^2,
\end{equation}
provided $\aaw\in(0,1)$ and $\lah_r\ge -\sfrac16\aaw$.
Consequently we find
\begin{equation}\lbl{i.mobd1}
|\mosym(\lah,\xi)| \le \frac{6}{\aaw|\xi|}, 
\qquad\mbox{$\xi\in\R+i\aaw$, \ $\re\lah\ge -\sfrac16\aaw$.}
\end{equation}
In particular, since 
$|\xi|\ge\aaw$,  
the left-hand side of \qref{i.mo1} is uniformly bounded 
away from zero, and $\mosym$ is uniformly bounded.

2. (Low frequencies) 
Now we carefully identify a long-wave regime where the 
result of Taylor expansion in \qref{e.expand3} 
yields the limits \qref{l.memo} and \qref{l.memoB}.
Fix $\nu_0\in(\frac13,\frac12)$ and let
\begin{equation}\lbl{d.I1}
I_0=\{\xi\in\R+i\aaw: |\eph\xi|\le \eph^{\nu_0}\}.
\end{equation}
Put $D_0=\mosym(\lah,\xi)\inv$, $E=\mesym(\lah,\xi)\inv-D_0$.
Then by \qref{e.expand3}, $E=\xi^2O(\eph^2\xi^2)$ and since $|D_0|\ge
\frac16\aaw|\xi|$ we find that uniformly for 
$\xi\in I_0$, $\lah\in\Omeh$ we have
\begin{equation}\lbl{i.mebd1}
 |\mosym(\lah,\xi)-\mesym(\lah,\xi)| 
=\left| \frac{E}{D_0(D_0+E)}\right| \le 
\left(\frac6{\aaw}\right)^2 
\frac{ K|\eph\xi|^2}{1-|\xi|K|\eph\xi|^2} \le \hat K \eph^{2\nu_0}
\end{equation}
for small enough $\eph$, since $\eph^2|\xi|^3\le \eph^{3\nu_0-1}=o(1)$.
Moreover, for $\xi\in I_0$ one has 
\[
\left|\sqrt{\frac{\gamma\tanh\eph\xi}{\eph\xi}} -1\right|
\le K\eph^{2\nu_0}.
\]
It follows that \qref{l.memoB} holds uniformly 
in this regime, as well.

3. For high frequencies the KdV limit is not relevant.
In this regime, the symbols $\mesym$ and $\mosym$ must be
shown separately to be small. Let us consider $\mosym$ first.
When $\xi\in I_0^c:=\R+i\aaw\setminus I_0$ we have $|\xi|\ge
\eph^{\nu_0-1}$, and 
from the estimate \qref{i.mobd1} it is clear that 
\begin{equation}\lbl{i.mobd2}
\sup_{\xi\in I_0^c,\,\lah\in\Omeh}
|\mosym(\lah,\xi)|\le \frac{6\eph^{1-\nu_0}}{\alpha} ,
\end{equation}
and this tends to zero as $\eps\to0$.

What remains to show is that 
\begin{equation}\lbl{i.mebd2}
\sup_{\xi\in I_0^c,\,\lah\in\Omeh}
 |\mesym(\lah,\xi)| 
\to 0 
\qquad\mbox{as $\eps\to0$.
}
\end{equation}
Since the square-root factor 
in \qref{l.memoB} is bounded, the proof of both
\qref{l.memo} and \qref{l.memoB} will be complete once
\qref{i.mebd2} is established.
This estimate is the most subtle of the symbol estimates.
It has nothing to do with the KdV limit, but rather expresses a uniform
stability property that holds at high frequency over moderately long 
time scales.  Its proof breaks into two further regimes for $|\xi|$ 
and involves using Corollary~\ref{c.DBsym} to interpolate between low
frequencies and high.

4. (High frequencies) Here we consider the set
\begin{equation}\lbl{d.i3}
I_\infty = \{\xi\in\R+i\aaw : |\eph\xi|\ge  K_2\}.
\end{equation}
where $K_2$ is large. 
In this regime, the denominator of $\mesym$ is estimated from below by
\[
|\eph\xi| -|\gamma\eph\xi\tanh\eph\xi|^{1/2} -|\eph^3\lah|
\ge |\eph\xi| -|\eph\xi|^{1/2} -1  
\ge \sfrac12|\eph\xi|,
\]
and consequently 
\begin{equation} \lbl{i.mebdc3}  
\sup_{\xi\in I_\infty,\,\lah\in\Omeh}
|\mesym(\lah,\xi)| 
\le 2\eph^2.
\end{equation}

5. (Transition frequencies)  
Fix $\nu_1$ with $\nu_0<\nu_1<\frac12$ and let
\begin{equation}\lbl{d.I2}
I_1 = \{\xi=k+i\aaw: 3\eph^{\nu_1}\le |\eph k|\le K_2\}.
\end{equation}
Recalling $\aw=\eph\aaw$ and $\gamma<1-\aw^2$, 
we apply Corollary~\ref{c.DBsym} together with the bound
 \qref{i.tanh2} valid for $|\kappa|$ small.
Then we find that $\eph$ small enough, 
with $\xi=k+i\aaw\in I_1$ and $-\lah_r<1$, 
the real part of the denominator of $\mesym(\lah,\xi)$
is negative and bounded (away from zero) by
\[
\re(-\eph^3\lah+i\eph\xi+\sqrt{-\gamma\eph\xi\tanh\eph\xi})
\le -\eph^3+\eph\aaw(-1+1-\eph^{2\nu_1})
\le 
-\sfrac12\eph^{1+2\nu_1}.
\]
By consequence we find that 
\begin{equation}\lbl{i.mebdc1}  
\sup_{\xi\in I_1,\,\lah\in\Omeh}
|\mesym(\lah,\xi)|  
\le 2K_2 \eph^{2-1-2\nu_1},
\end{equation}
and this tends to zero since $\nu_1<\frac12$.

%
Now $\R+i\aaw=I_0\cup I_1\cup I_\infty$, and 
the outstanding estimate \qref{i.mebd2} is established.
This finishes the proof of the limits 
\qref{l.memo}-\qref{l.memoB}.

{\bf Remark.} The analysis of symbol limits in this section is 
simpler than the one carried out for lattice solitary waves in \cite{FP4}.
Partly this is due to the simple way that $\lah$ appears
in the denominator of $\mesym$ here. 
But partly it is due to the fact that here we must study the regime
$|\eph^3\lah|\ge \hat K$ by other means, since
the symbol $\mesym$ in \qref{d.mesym} is {\it not} bounded on 
the whole set where $\re\lah\ge0$ and $\xi\in\R+i\aaw$.

\subsection{Limit of junk terms}\lbl{s.jlim}
To complete the proof of Theorem~\ref{t.wlim}
it suffices to prove \qref{l.w3}, i.e., show that $\|J_*\|_a=o(1)$ 
in operator norm on $L^2_a$,
uniformly for $\lambda\in\Omee$ where
\begin{equation}\lbl{c.lam2}
\Omee=\eph^3\Omeh= \{\lam\in\C: |\lam|\le 1,\ \re\lam\ge -\sfrac16\aaw\eph^3\}.
\end{equation}

By the bound
\qref{i.j1221} and the resolvent estimate for $\AAp$ in 
Corollary~\ref{c.DBsym}, we have
\begin{equation}\lbl{i.Js2}
\|(\lam-\AAp)\inv J_{12}(\lam-\Ac_{22})\inv J_{21}\|_a \le K \eph^2.
\end{equation}
Further, Corollaries~\ref{c.DBsym} and \ref{c.comm1} imply 
$\|(\lam-\AAp)\inv J_{11}\|_a\le K$, hence
\begin{equation}\lbl{i.Jst}
\|J_*\|_a\le K 
\end{equation}
uniformly for $\re\lam\ge-\frac16\aaw\eph^3$.
It remains to prove that uniformly for $\lambda\in\Omee$,
\begin{equation}\lbl{i.j11a}
\|(\lam-\AAp)\inv J_{11}\|_a \to 0 \quad\mbox{as $\eph\to0$}.
\end{equation}
Then $\|J_*\|_a=o(1)$ will follow.
(We remark that it appears this may not be true uniformly 
for all $\lam$ satisfying $\re\lam\ge0$, however.)

Because of the estimates \qref{i.Js2} and \qref{i.AAp}
and the expression for $J_{11}$ in \qref{d.junkmat}, 
it suffices to prove that the operators
\begin{eqnarray}
J_0 &=& (\lam-\AAp)\inv(u_q'\sqrt\gamma +[\Sc,u_q]),\\
J_1 &=& (\lam-\AAp)\inv(u_p'-[\Sc,u_p]\SC\inv\D), 
\end{eqnarray}
are $o(1)$ in $L^2_a$ operator norm 
uniformly for $\lam\in\Omee$ --- 
note that $|q'pq\inv-\sqrt\gamma u_q'|=O(\eph^5)$.

In order to bound $J_0$ it would suffice to note $u_q'=[\D,u_q]$
and apply Proposition~\ref{P.comm} with $g=u_q$ and with appropriate
symbols.
However, the form of $J_1$ is slightly different, so what we do
instead is observe that the weight-transformed
operators $\tilde J_j=\emax J_j\eax$ ($j=0, 1$)
act on a given smooth $f$ with compact support via
\begin{equation}
(\FF \tilde J_j f)(k) = \sqrt{\gamma^{1-j}}
\intR P(k)\Big(
i(k-\hat k)+(Q(k)-Q(\hat k))R(\hat k)^j 
\Big)
\hat g(k-\hat k) \hat f(\hat k)\,\frac{d\hat k}{2\pi}
\end{equation}
for $j=0$ and 1, with
\[
P(k) = \frac1{-\lam+i\xi+\sqrt{-\gamma\xi\tanh\xi}},
\qquad
Q(k) = \sqrt{-\xi\tanh\xi},
\qquad
R(k)= \sqrt{\frac{\xi}{\tanh\xi}}.
\]
Here $\xi=k+ia$, and $g=u_p$ or $u_q$ has the form $g(x)=\eph^2G(\eph x)$
with $G$ bounded in $H^2$ as in section~\ref{S.comm}.
By almost the same short proof as that of Proposition~\ref{P.comm},
we find that 
\begin{equation}\lbl{i.Jj}
\|J_j f\|_a \le C_* C_G \|f\|_a
\end{equation}
where $C_G$ is as in Proposition~\ref{P.comm} (and is uniformly bounded), and
\begin{eqnarray}
C_* &=& \sup_{k,\kh\in\R} \eph^2 
\frac{|P(\eph k)||i\eph(k-\kh)+(Q(\eph k)-Q(\eph\kh))R(\eph\kh)^j|}
{\ip{k-\kh}^s}
\nonumber\\
&=& \sup_{k,\kh\in\R} \eph^3|P(\eph k)| 
\left|1+ 
\frac{Q(\eph k)-Q(\eph\kh)}{i(\eph k-\eph\kh)} R(\eph\kh)^j\right|
\frac{|k-\kh|}{\ip{k-\kh}^s}.
\lbl{d.cs2}
\end{eqnarray}
Here we take $s=\frac43$ as previously. This implies that the 
last factor in \qref{d.cs2} is bounded.
To bound the other factors we consider the case $|\eph\xi|\le \eph^\nuh$
and its opposite, for any $\nuh\in(0,1)$ fixed.

In the first case, $|\eph\xi|\le\eph^\nuh$, the factor 
$|P(\eph k)|\eph^3\le K$ uniformly since by Corollary~\ref{c.DBsym},
the real part of the denominator of $P(\eph k)$
is bounded away from zero, 
satisfying $\re P(\eph k)\inv\le-\frac14\eph^3\aaw$.
The middle factor is $O(\eph^\nuh)$ and tends to zero uniformly, since
the symbols $Q$ and $R$ are analytic near $\xi=0$ and
$Q'(k)\to -i$ and $R(k)\to1$ as $\xi=k+ia\to0$.

In the other case, $|\eph\xi|\ge\eph^\nuh$, we note that the middle factor
is uniformly bounded due to the estimate \qref{i.Rlip}.
Due to Corollary~\ref{c.DBsym}, 
\[
\re P(\eph k)\inv \le
-\re\lam+\eph\aaw\left(-1+ \sqrt{\frac{\tanh \eph^\nuh}{\eph^\nuh}}
\right) \le 
\aaw\eph^3-\sfrac19{\aaw \eph^{1+2\nuh}}.
\]
Hence $|P(\eph k)|\eph^3\le K\eph^{2-2\nuh}$ for small $\eph$, 
and we conclude that 
\begin{equation}
C_* \le K(\eph^{1+2\nuh}+\eph^\nuh)
\end{equation}
which tends to zero uniformly for $\lam\in\Omee$. 

\section{Analysis of the bundle limit}\lbl{S.gsr}

For the remainder of the proof of our asymptotic linear stability theorem,
there are two approaches possible. One is to proceed in a fashion
similar to the treatment of FPU lattice waves in the KdV limit in
\cite{FP4}. In that approach, one notes that any eigenfunction of $\Ac$
corresponding to a nonzero eigenvalue is orthogonal to two particular
elements of the generalized kernel of the adjoint $\Ac^*$.  (In \cite{FP4}
this was expressed in terms of symplectic orthogonality, using
Hamiltonian structure.) This yields reduced orthogonality conditions
that are necessary for elements of the kernel of the scalar bundle
$W(\lam)$. After an appropriate scaling, one proves convergence of
these conditions to corresponding ones for the KdV bundle $W_0(\lah)$,
in a dual space.  Then uniform invertibility of $W(\lam)$ on the
codimension-2 subspace satisfying the orthogonality conditions follows
by a straightforward perturbation argument.

We prefer to emphasize, however, that the required spectral properties follow
from the bundle convergence theorem~\ref{t.wlim}  by `soft' arguments
based on the GSR perturbation theorem,
and does not require further convergence analysis of adjoint zero
modes. There are essentially only two `hard' points left. 
Namely, we need to show that (i) the bundle $W(\lam)$ is Fredholm of index
zero for relevant values of $\lambda$, and 
(ii) the solitary-wave degrees of freedom (translational shift and wave speed)
naturally provide two independent elements in the generalized 
kernel of $\Ac$.  
In comparing the need for point (i) with the alternative approach, 
we observe that if one knows $\lam-\Ac$ has empty kernel,
one would likely need to prove a Fredholm property anyway
to conclude that $\lam-\Ac$ is surjective and $\lam$ is in
the resolvent set of $\Ac$.

In this section we will establish point (i), and invoke 
Gohberg-Sigal-Rouch\'e perturbation theory 
to characterize the null multiplicity of characteristic 
values of the bundle $W(\lam)$. 
This is related to the algebraic multiplicity of eigenvalues of $\Ac$
in the following section.
Point (ii) is dealt with in Appendix B.

\subsection{Fredholm property of the bundle}

\begin{lemma}\lbl{L.Fred}
For $\eph>0$ sufficiently small, 
$W(\lam)$ is Fredholm with index 0 for all $\lam\in\Omee$. 
\end{lemma}

\noindent{\it Proof.} It suffices to show that we may write
\begin{equation}\lbl{e.Weq1}
W(\lam)=W_i+W_c
\end{equation}
where $W_i$ is invertible and $W_c$ is compact. 
To demonstrate this, we use Fourier filters $\pi_o$ and $\pi_i$
as defined in \qref{d.0lohi}-\qref{d.lohi}. It is convenient however
to use a soft wavenumber cutoff in the range $[\kah,2\kah]$ 
with $\kah=\eph^\nu$ where $0<\nu<\frac35$  (see \qref{d.I1}).
More precisely, fix $\phi(k)=1$ ($|k|\le1$), $2-|k|$ ($1\le|k|\le2$),
$0$ ($|k|\ge2$) and set $\phi_\eps(k)=\phi(k/\eph^\nu)$ and 
in place of \qref{d.0lohi}-\qref{d.lohi} define
\begin{equation}\lbl{d.lohi2}
\pi_{0o}= \FF\inv \phi_\eps\FF 
\qquad
\pi_o = \emax\pi_{0o}\eax, 
\qquad
\pi_i=I-\pi_o.
\end{equation}
We then define
\begin{equation}\lbl{d.Wc}
W_i=I-\pi_i(\lam-\AAp)\inv\UU-J_*,
\qquad
W_c = -\pi_o(\lam-\AAp)\inv\UU.
\end{equation}
The operator $W_i$ is uniformly invertible for $\lam\in\Omee$.
This is so since $\|J_*\|_a$ is uniformly small and so is the
middle term, for we have $\|\D\inv\UU\|_a\le K\eph^2$, while
$\emax\pi_i(\lam-\AAp)\inv\D\eax$ is a Fourier multiplier with symbol
dominated by $\mesym$, with 
\begin{equation}\lbl{i.fred}
\|\eph^2\pi_i(\lam-\AAp)\inv\D\|_a\le 
\sup_{\xi\in I_0^c,\,\lah\in\Omega_\eps}
 |\mesym(\lah,\xi)|  \to 0
\end{equation}
as $\eps\to0$ due to \qref{i.mebd2}.
Then, if $\eps$ is small enough, $W_i$ is invertible 
for all $\lam\in\Omee$.

On the other hand, the operator $W_c$ on $L^2_a$ 
is equivalent to the weight-transformed $\emax W_c\eax$ on $L^2$.
The latter operator is compact
by the convenient compactness criterion of \cite{Pe85}---It
is the sum of two terms
of the form $\FF\inv\phi_1\FF\phi_2$, where $\phi_1$ and $\phi_2$
are multipliers by bounded continuous functions on $\R$ that 
approach zero at infinity. This finishes the proof of the Lemma.

{\bf Remark.} We note that in the decomposition \qref{e.Weq1},
both terms $W_c$ and $W_i$ are analytic functions of $\lam$ for
$\lam\in\Omee$.  (This fact will be used in studying the
full resolvent of $\Ac$.)

\subsection{Characteristic values and the Gohberg-Sigal-Rouch\'e theorem}
We first recall some relevant basic information from \cite{GS71}.
(We change some terminology slightly for clarity. 
An alternative source is \cite{GGK90}.)
Let $X$ be a Hilbert space, and suppose a function $\lam\mapsto\WW(\lam)$ 
is analytic on a complex domain $\Omega_0\subset\C$, taking values
in the space of bounded linear operators on $X$, and all its values
are Fredholm of index zero. A point $\lam_0$ is a {\it characteristic
value} of $\WW$ if $\WW(\lam_0)$ has a nontrivial kernel.
A {\it root vector} is an analytic function $z(\lam)$ with values in $X$
satisfying $\WW(\lam_0)z(\lam_0)=0$ with $z(\lam_0)\ne0$. 
The {\it order} of a root vector at $\lam_0$
is the order of $\lam_0$ as a zero of $\WW(\lam)z(\lam)$.
The {\it null multiplicity} of a characteristic value is a positive
integer whose precise definition in general need not concern us here.
The {null multiplicity} of $\lam_0$ is always at least as large as the 
maximum order of any root vector.
Furthermore, the null multiplicity equals this maximum order if
and only if the kernel of $\WW(\lam_0)$ is one-dimensional. 

Suppose $\Omega$ is a subdomain of $\Omega_0$, with boundary
$\Gamma$ that is a simple closed rectifiable contour in $\Omega_0$,
and suppose $\WW(\lam)$ is invertible for all $\lam\in\Gamma$.  
The sum of all null multiplicities for all characteristic values in 
$\Omega$ is denoted $n(\WW,\Omega)$ and is called the 
{\it total multiplicity of $\WW$ in $\Omega$}.  
A simple corollary of a far-reaching generalization of 
Rouch\'e's theorem proved by Gohberg and Sigal \cite{GS71,GGK90} is the following.
\begin{theorem} Assume that for $j=1$ and 2, $\WW_j(\lam)$ is
analytic and Fredholm of index zero in $\Omega\cup\Gamma$. 
Assume that for all $\lam\in\Gamma$,
$\WW_1(\lam)$ is invertible and the operator norm
\[
\|\WW_1(\lam)\inv(\WW_1(\lam)-\WW_2(\lam)\|_{X} <1.
\]
Then $\WW_2(\lam)$ is invertible on $\Gamma$, and the total
multiplicity $n(\WW_2,\Omega)=n(\WW_1,\Omega)$.
\end{theorem}

We apply this abstract result with $\WW_1=W_0$ as defined in \qref{d.W0},
and $\WW_2=\tiW$ as defined in Proposition~\ref{t.wlim}.
We take $X=(L^2_a)^2$, and the contour $\Gamma$ as the boundary of the set
$\Omega=\Omeh$.
As a consequence of Proposition~\ref{t.wlim}, for $\eph>0$ sufficiently
small, the null multiplicity $n(\tiW,\Omeh)=n(W_0,\Omeh)$.
But the latter number is 2, as a consequence of the following result.
\begin{prop}\lbl{p.kdv} Suppose $0<\aaw<\sqrt3$ and 
$\beta=\frac12\aaw(1-\frac13\aaw^2)$.  In $L^2_\aaw$,
the only characteristic value of $W_0(\lah)$ with
$\re\lah>-\beta$ is $\lah=0$, and this value has null multiplicity 2.
\end{prop}
This Proposition mainly follows from known facts concerning 
the eigenvalue problem for the KdV soliton. We provide a self-contained 
proof in Appendix C for the reader's convenience.

\begin{cor}\lbl{c.nmul}
For $\eph>0$ sufficiently small, $W(\lam)$ is invertible
for all $\lam$ on the boundary of $\Omee=\eph^3\Omeh$, and the 
total multiplicity of $W$ in $\Omee$ is 2.
\end{cor}


\section{Analysis of resolvent and eigenvalues}\lbl{S.res}

It remains to complete the proof that for $\eph>0$ sufficiently
small, the operator $\Ac$ has no eigenvalue with
$\re\lam\ge-\frac16\aaw\eph^3$ other than $\lam=0$, which
is a discrete eigenvalue with algebraic multiplicity 2.
Conditional asymptotic stability will then follow directly 
from the Gearhart-Pr\"uss theorem.

\subsection{Resolvent and spectral projection}

To begin, we note that by simple elimination, 
whenever $\lam\in\Omee$, the resolvent equation
\begin{equation}\lbl{e.res0}
(\lam-\Ac)
\pmat{f_1 \\ f_2} 
=
\pmat{g_1 \\ g_2}
\end{equation}
is equivalent to 
\begin{eqnarray}
\lbl{e.res.1}
W(\lam)f_1 &=& (\lam-\Ac_+)\inv(g_1+J_{12}(\lam-\Ac_{22})\inv g_2), \\
\lbl{e.res.2}
f_2 &=& (\lam-\Ac_{22})\inv (J_{21}f_1+g_2).
\end{eqnarray}
Then clearly, $\lam$ is in the resolvent set of $\Ac$
if $W(\lam)$ is invertible, so any point of the spectrum of $\Ac$
inside $\Omee$ must be a characteristic value of $W(\lam)$.
(The converse is not so easy to argue, since a root vector in $L^2_a$ 
need not be in $H^1_a$. We finesse this point in the following argument.)

For $\eph>0$ small, as a consequence of Corollary~\ref{c.nmul},
there are at most 2 points of $\Omee$ in the spectrum of $\Ac$. 
We claim that each such spectral point $\lam_0$ is a {\it discrete
eigenvalue} of $\Ac$, which means that the associated spectral projection,  
\begin{equation}\lbl{d.P0}
P_0 = \frac{1}{2\pi i}\int_{\Gamma_0} (\lam-\Ac)\inv\,d\lam,
\end{equation}
has finite rank.
Here $\Gamma_0$ is a small enough circle about $\lam_0$
enclosing no other point of the spectrum.
To prove this, note that from the decomposition
formula \qref{e.Weq1}, we can write
\begin{equation}\lbl{e.resA1}
W(\lam)\inv = W_i\inv - W(\lam)\inv W_c W_i\inv,
\end{equation}
from which we easily deduce from \qref{e.res.1}-\qref{e.res.2} that 
we can write $(\lam-\Ac)\inv = \RR_i+\RR_c$
where $\RR_i$ is analytic in $\Omee$ and $\RR_c$ is compact.
Then  the integral $\int_{\Gamma_0}\RR_i\,d\lam=0$ and it
follows that $P_0$ is compact. Since $P_0$ is a projection,
it has finite rank, and it follows that its range consists entirely of
generalized eigenvectors of $\Ac$.

\subsection{Algebraic multiplicity of eigenvalues}
It remains to relate the algebraic
multiplicity of an eigenvalue $\lam_0$ of $\Ac$
to the null multiplicity of $\lam_0$ as a characteristic value
of $W(\lam)$.  These quantities are in fact equal, 
but for present purposes it suffices to be brief and
prove a simpler, weaker result.

\begin{prop} For $\eps>0$ sufficiently small,
if $\lam_0\in\Omee$ is an eigenvalue of $\Ac$, then 
$\lam_0$ is a characteristic value of $W$.
Furthermore, if a Jordan chain $z_1,\ldots,z_k$ is a
Jordan chain of elements in $(H^1_a)^2$ satisfying
\[
(\Ac-\lam_0)z_j=z_{j-1} 
\quad\mbox{ for $j=1,\ldots,k$, with $z_0=0$,}  
\]
then a root vector $\eta(\lam)$ of order at least $k$ exists
for $W$ at $\lam_0$.
\end{prop}

\noindent {\it Proof.} Supposing $z_1,\ldots,z_k$ 
is a Jordan chain for $\Ac$ of length $k$, let 
$f(\lam)=\sum_{j=1}^k  (\lam-\lam_0)^{j-1}z_j$.
Then $f(\lam)$ is analytic with values in $(H^1_a)^2$ (the domain of
$\Ac$) and $(\lam-\Ac)f(\lam)=(\lam-\lam_0)^kz_k=:g(\lam)$.
By elimination, \qref{e.res.1}-\qref{e.res.2} hold, and
consequently $W(\lam)f_1(\lam)=O(|\lam-\lam_0|^k)$. 
Thus there is a root vector $\eta(\lam)=f_1(\lam)$ of order 
at least $k$, and $\lam_0$ is a characteristic value of $W$.

\subsection{Proof of asymptotic stability}
Recall $\aaw\in(0,\half]$ is fixed, and take $\eph>0$ sufficiently small.
As a consequence of the last Proposition and the fact from Appendix B
that $\lam=0$ has algebraic multiplicity {\it at least} 2 
for $\Ac$, we conclude that the null multiplicity of the characteristic value 
$\lam=0$ for the bundle $W$ is at least 2.  Since the total multiplicity of
the bundle $W$ in $\Omee$ is 2, we deduce that
(i) there are no nonzero points of the spectrum of $\Ac$ in $\Omee$,
and (ii) the algebraic multiplicity of $\lam=0$
is exactly 2. In particular, the kernel of $\Ac$ is simple
(and the same is true for $W(0)$).

Consequently, the spectral projection $P_0$ for $\lam=0$ has rank 2,
and restricted to the complementary invariant subspace
$\bar Y_a=(I-P_0)(L^2_a)^2=\ker P_0$, 
the resolvent $(\lam-\Ac)\inv$ is bounded
uniformly for $\lam\in\Omee$. By consequence of
Proposition~\ref{p.resA}, this restricted resolvent is bounded
uniformly for all $\lam\in\C$ with $\re\lam\ge -\sfrac16\aaw\eph^3$.
It follows automatically that the restricted resolvent is bounded uniformly 
in a slightly larger half-plane $\re\lam\ge -\beta$ for some 
$\beta>\sfrac16\aaw\eph^3$.
Using the Gearhart-Pr\"uss asymptotic stability criterion
(see Corollary 4 in \cite{Pr85})
gives us the conditional linear asymptotic stability result 
claimed in Theorem~\ref{t.Amain}.

\section{Spectral stability without weight}

In this section we prove Theorem~\ref{t.nowt2}, showing that in the 
unweighted space $(L^2)^2$, 
the spectrum of the operator $\Ac$ is the imaginary axis,
if $\eph>0$ is sufficiently small.
The proof breaks into four steps. For $\re\lam>0$, we show that (i) either 
$\lam$ is in the resolvent set or $\lam$ is an eigenvalue, and (ii) if $\lam$ is
an eigenvalue in $(L^2)^2$, then it is an eigenvalue in $(L^2_a)^2$.
Since by Theorem~\ref{t.main} there are no such eigenvalues,
this proves that $\Ac$ has no spectrum in the right half-plane.
Next we show that (iii) $\Ac$ has no spectrum in the left half-plane
due to a symmetry under space and time reversal. Finally, we show
(iv) each point of the imaginary axis does belong to the spectrum of
$\Ac$, by a fairly standard construction of a sequence of approximate
eigenfunctions.

1. Suppose $\re\lam>0$. To accomplish the first step, as in section 9 we write 
\[
\Ac=\Ac_*+\RR_*, \qquad
\Ac_*= \pmat{\AAsp &0\\0& \AAsm}, \qquad
\RR_* = \pmat{\ti J_{11} & J_{12}\\ J_{21} & J_{22}}.
\]
By applying Lemma~\ref{l.AAspm} with $\aaw=0$,  we infer
that $\lam$ belongs to the resolvent set of $\Ac_*$ and 
\[
\lam-\Ac = 
(I-\RR_*(\lam-\Ac_*)\inv) (\lam-\Ac_*)
\]
We claim that $\RR_*(\lam-\Ac_*)\inv$ is compact, whence it follows that
either $\lam$ is in the resolvent set of $\Ac$ or it is an eigenvalue.
To prove the claim, it suffices to show that each entry is a sum of 
terms each of which is a product of bounded operators,
at least one of which is compact.
Let $L=I+\D$ and note that since the domain of $\AAspm$ is $H^1$,
the operators $L(\lam-\AAspm)\inv$ are bounded on $L^2$.
Thus it suffices to show that $\RR_*L\inv$ is compact.
By the criterion of \cite{Pe85}, an operator of the form
$g \QQ$ or $\QQ g$ is compact on $L^2$ provided that $g$ is a pointwise multiplier
by a continuous function satisfying $g(x)\to0$ as $|x|\to\infty$,
and $\QQ$ is a Fourier multiplier with continuous symbol satisfying
$\hat\QQ(k)\to 0$ as $|k|\to\infty$.

We now deal with the various terms in $\RR_*L\inv$ 
from \qref{d.junkmat}-\qref{d.A11b}.
With the notation of section 6, note $q'=u_q'$ decays as $|x|\to\infty$,
and $L\inv$ has symbol $(1+ik)\inv$ tending to 0 as $|k|\to\infty$.
Hence the operator $R_1L\inv=pq\inv q'L\inv$ is compact. 
Similarly $p'L\inv$ is compact. 

To treat terms involving commutators, we consider first
the worst term, $[p,\SC] \SC\inv\D L\inv$.
Note that the operator $\SC\inv\D L\inv$ is a Fourier multiplier
with bounded continuous symbol $\sqrt{k/\tanh k}/(1+ik)$, 
hence is bounded.
(The symbol decays, but we do not use this fact.)
We now claim that 
\begin{equation}
[p,\SC] \mbox{\ is compact.}
\end{equation}  
We will show, in fact, that $[p,\SC]$
is the uniform-norm limit of a sequence $[p,\SC_n]$, where $\SC_n$ is a 
Fourier multiplier with continuous symbol of compact support.
Let $\phi:\R\to[0,1]$ be a smooth cutoff function, taking the value 1 
on $[-1,1]$ and 0 on $\R\setminus[-2,2]$ and let $\psi=1-\phi$. 
Let $\SC_n$ be the Fourier multiplier with symbol 
$\phi(k/n)\sqrt{-\tig k\tanh k}$. 
Then $u_p\SC_n$ and $\SC_nu_p$ are each compact.
As in the proof of Corollary~\ref{c.comm1},
using Proposition~\ref{P.comm}, the $L^2$ operator norm
of $[p,\SC]-[p,\SC_n]=[u_p,\SC-\SC_n]$ is bounded by $\eph^3C_nC_G$, where 
$C_G$ is bounded and 
\[
C_n = \sup_{k,\kh\in\R} 
\frac{|Q(k)\psi(k/n)-Q(\kh)\psi(\kh/n)|}{\ip{k-\kh}^{4/3} },
\qquad  Q(k)=\sqrt{k\tanh k}
\]
Since $Q$ is increasing and $Q'$ decreasing for large $k$, we have
the uniform derivative estimate
\[
|(Q(k)\psi(k/n))'|= |Q'(k)\psi(k/n)+Q(k)\psi'(k/n)/n|\le |Q'(n)|+K|Q(2n)/n|\le K/\sqrt{n}.
\]
Then it easily follows $C_n\le K/\sqrt n\to0$ as $n\to\infty$.
Hence $[p,\SC]$ is compact on $L^2$.

Similarly the commutators $[\SC,q]$ and $[\SC,\rho]$ are compact,
and it follows directly that $\RR_* L\inv$ is compact.
This finishes the first step.

2. For the second step, suppose $\lam$ is an eigenvalue with $\re\lam>0$
and with eigenfunction $(\eta_4,\phi_4)\in (H^1)^2$, the domain of $\Ac$.
We then can write (from \qref{e.ev4b})
\[
\pmat{\lam-\AAp&0\\0&\lam-\AAm}
\pmat{\eta_4\\ \phi_4} 
= \pmat{g_1\\ g_2},
\]
with
\[
\pmat{g_1\\ g_2}= 
-\frac12\pmat{
-\D u_p
-\SC u_q-u_q\SC +R_1+R_2 &
-\SC u_q+u_q\SC +R_1-R_2 \cr
\SC u_q-u_q\SC +R_1-R_2 &
-\D u_p+
\SC u_q+u_q\SC +R_1+R_2 
}
\pmat{\eta_4\\ \phi_4} .
\]
We claim that $g_1$ and $g_2$ lie in $L^2_a$ as well as $L^2$.
This is not difficult to check, since $\emax u_p$ and $\emax u_q$ are in $H^2$. 
Since $\re\lam>0$, the Fourier multipliers $(\lam-\AApm)\inv$
are bounded on $L^2$ and on $L^2_a$. Indeed, they map the subspace
$L^2\cap L^2_a$ of $L^2$ into $H^1\cap H^1_a$ 
(as one can check by approximation using smooth test functions
and analyticity of the Fourier transform for $0<\im \xi< a$). 
It follows that $(\eta_4,\phi_4)\in (H^1_a)^2$, and that
$\lam$ is an eigenvalue of $\Ac$ in the space $(L^2_a)^2$.
But there is no such eigenvalue for $\eph>0$ sufficiently small,
by Theorem~\ref{t.main}. This concludes the proof of spectral stability
for $\Ac$ in $(L^2)^2$.

3. The resolvent equation for $\Ac$ has a symmetry under space and
time reversal inherited from the original water wave equations.
For present purposes, 
this is most easily studied in terms of the variables used in
\qref{e.ev3},
for which the resolvent equation may be written (in $L^2\times L^2$)
\begin{equation}\lbl{e.rs3}
\pmat{\lam-q\D pq\inv & q\SC\\ \SC q &\lam-\SC p\D\SC\inv}
\pmat{\eta_3\\ \phi_3}=
\pmat{f_1\\ f_2}.
\end{equation} 
Recall that $p$ and $q$ are even functions. 
Let $\CC$ be the space reversal operator, $\CC f(x)=f(-x)$.
Now, $\SC$ preserves parity (since its symbol is even)
while $\D$ reverses it.
Applying space reversal to \qref{e.rs3} the problem is seen to be
of the same type after the replacements
\[
\lam\mapsto -\lam, \qquad
\pmat{\eta_3\\ \phi_3}\mapsto \pmat{-\CC \eta_3\\ \CC \phi_3}.
\]
Thus $\lam$ is in the resolvent set if and only if $-\lam$ is.
It follows at this point that the $L^2$ spectrum of $\Ac$ is contained 
in the imaginary axis.

4. Suppose $\re\lam=0$. Then there exists $\kh\in\R$ such
that $\lam=\AAp(\kh)=i\kh+i\sqrt{\tig \kh\tanh\kh}$.
Formally, $(\lam-\AAp)e^{i\kh x}=0$. 
We construct a sequence of approximate eigenfunctions for $\Ac$ in 
$(L^2)^2$ by a cutoff and translation argument.
Fix a smooth function $\psi$ with compact support, and consider
pairs $(\eta_4,\phi_4)$ of the form
\[
\eta_4= e^{i\kh(x+\tau)}\psi(\eh(x+\tau))\sqrt{\eh},
\qquad \phi_4=0.
\]
The $L^2$ norm of $\eta_4$ is independent of $\eh$ and $\tau$.
We claim that taking $\eh=1/n$, we can choose $\tau$ depending on
$n$ such that 
$\|(\lam-\Ac)(\eta_4,0)\|_{L^2}\to0$ as $n\to\infty$.
Due to the structure of $\Ac$ in \qref{e.evfinal} it suffices to show
that as $n\to\infty$. in $L^2$ we have
(a) $(\lam-\AAp)\eta_4\to0$, and (b) 
$(\D u_p +\SC u_q+J_{11})\eta_4$ and $J_{21}\eta_4 \to0$.

To prove (a), we simply note that the Fourier transform
\[
\FF((\lam-\AAp)\eta_4) (k) = e^{ik\tau}(\AAp(\kh)-\AAp(k))
\hat\psi\left(\frac{k-\kh}{\eh}\right)\frac1{\sqrt{\eh}},
\]
and this tends to 0 in $L^2$ as $\eh\to0$ uniformly in $\tau$.

To prove (b), it is convenient to note that for any fixed $\eh>0$,
$\eax\eta_4\to0$ in $H^2$ as $\tau\to\infty$.
Moreover, $u_p\emax$ and $u_q\emax$ are bounded in $H^1$.
Then it follows, for example, that $u_p\eta_4$ and $u_q\eta_4\to0$
in $H^1$ as $\tau\to\infty$, and 
\[
R_2\eta_4 = \Bigl(\D (u_p\emax)- \SC(u_p\emax)
(\eax\D\SC L\inv\emax)
\Bigr) (\eax\eta_4) \to 0
\]
in $L^2$ as $\tau\to\infty$, since the weight-transformed 
operator $\eax\D\SC L\inv\emax$ has bounded symbol
and is bounded on $H^1$. ($L=1+\D$ as above.)
Similarly it follows $[S,u_q]\eta_4$ and $R_1\eta_4\to0$
in $L^2$ as $\tau\to\infty$. Choosing $\tau$ appropriately
depending on $\eh$, this finishes the proof of (b).
Thus each point of the imaginary axis belongs to the $L^2$
spectrum of $\Ac$.

\appendix

\section{Rigorous asymptotics for solitary wave profiles}

Here we provide simple proofs of the estimates on the solitary wave
profile needed for our analysis of the eigenvalue problem and
resolvent.  For a sharper treatment of solitary water waves in
the limit $\eph\to0$ see Beale's work \cite{Be77}.

We work with the scaled form of \qref{e.uufp}, 
written in terms of $\us$ defined through
\begin{equation}\lbl{d.usa}
\uu(\xf)=\eph^2\us(\eph\xf).
\end{equation}
In terms of $\us$, we can write the fixed point equation \qref{e.uufp}
in the form
\begin{equation}\lbl{e.usfp}
\us = F(\us)=Q N(\us),
\end{equation}
with Fourier multiplier $Q$ and nonlinearity $N$ defined by
\begin{equation}
Q = \eph^2
\left(1-\gamma\frac{\tanh \eph D}{\eph D}\right)\inv,
\quad
N(\us) = \frac{\displaystyle 
\frac32\us^2 + \eph\us^3 + \frac12 (\tanh \eph D\,\us)^2
\left(1-2\gamma\eph \frac{\tanh \eph D}{\eph D}\us \right) 
}
{(1+\eph\us)^{2}}.
\end{equation}
We study this equation in a weighted Sobolev space of even functions.
For $\aaw>0$ fixed, let
\begin{equation}\lbl{d.ussp}
\Xm = \{f\colon\R\to\R: e^{\aaw x}f\in H^m,\  \mbox{$f$ even} \},
\qquad \Ym = \{f\colon\R\to\R: e^{\aaw x}f\in H^m,\  \mbox{$f$ odd} \},
\end{equation}
with the same norm (recall $\ip{k}=(1+k^2)^{1/2}$)
\[
\|f\|_{\Xm} = \|f\|_\Ym= \|e^{\aaw x}f\|_{H^m} = 
\left(\frac1{2\pi}\intR 
|\ip{k}^m\hat f(k+i\aaw)|^2\,dk\right)^{1/2}.
\]
One has $f\in\Xm$ (resp. $\Ym$) if and only if $f$ is even (resp. odd) and 
$\cosh\aaw x\, \D^jf\in L^2(\R)$ for $j=0,...,m$.
For $m\ge1$, the space $\Xm$ is a Banach algebra, while the bilinear product
map $(f,g)\mapsto fg$ is continuous from $\Ym\times\Ym$ to $\Xm$.
The intersection of exponentially weighted $H^m$ spaces is the direct sum of $\Xm$ and $\Ym$:
\[
H^m_\aaw\cap H^m_{-\aaw}
= \Xm\oplus\Ym.
\]

Due to the Taylor expansion of $\tanh$,
the symbol of $Q$ has the expansion
\begin{equation}\lbl{e.Sexp}
\hat Q(\xi)= \frac1{1+\frac13\gamma\xi^2+\xi^2O(\eph^2\xi^2)}, 
\qquad\xi=k+i\aaw.
\end{equation}
Formally, the limit of the fixed point equation \qref{e.usfp} is 
\begin{equation}\lbl{e.us0fp}
\us = Q_0 N_0(\us),
\qquad Q_0 = (1-\sfrac13\D^2)\inv, \qquad N_0(\us) = \sfrac32\us^2.
\end{equation}
Provided $0<\aaw<\sqrt3$, 
this fixed-point equation is satisfied by the KdV traveling-wave
profile
\begin{equation}\lbl{d.Sa}
\sw(x) = \sech^2(\sqrt3 x/2).
\end{equation}
This fixed point is nondegenerate in the space $\Xm$.
Indeed, the linearized map $\us\mapsto \us-Q_0(3\sw\us)$
has bounded inverse on $\Xm$, for the following reason.  
It is straightforward to show that the
map $Q_0\sw$ is compact on $H^m_\aaw\cap H^m_{-\aaw}$ (using \cite{Pe85}). So
if $\us\mapsto \us-Q_0(3\sw\us)$ is not an isomorphism on $\Xm$, then 
it vanishes for some nontrivial $\us$. By a simple bootstrapping argument,
this $\us$ must be a smooth function satisfying
$(I-\frac13\D^2+3\sw)\us=0$  with $e^{\aaw x}\us\in L^2$.
But only constant multiples of the odd function $\us=\sw'$ have this property;
from standard results for asymptotic behavior in ordinary differential equations,
any independent solution $\us$ grows like $e^{\sqrt3|x|}$ as $x\to\pm\infty$.

\begin{theorem}\lbl{t.swa}
Fix $m\ge2$, $\aaw\in(0,\sqrt3)$, and $\nus\in(0,1)$. Then
for $\eph>0$ sufficiently small, equation \qref{e.usfp} has 
a unique fixed point that belongs to $\Xm$ and satisfies
$\|\us-\sw\|_\Xm < \eph^\nus$.
This fixed point $\us$ depends smoothly on $\eph$.
\end{theorem}

To prove this result, we will invoke a standard fixed-point lemma
in the simple quantitative form from the appendix of \cite{FP1}. 
To make the estimates needed, we 
single out one difficult nonlinear term and let
\begin{equation}\lbl{d.N1}
N_1(\us)= (\tanh \eph D\,\us)^2. 
\end{equation}
Note $\tanh\eph D\,\us$ is odd if $\us$ is even, 
and $\tanh\eph D$ maps $\Xm$ to $\Ym$ continuously.  
Then we write
\begin{equation}\lbl{d.SNF}
Q=Q_0+Q_1, \qquad N=N_0+N_1+N_2, \qquad F=F_0+F_1+F_2+F_3, 
\end{equation}
with 
\[
F_0=Q_0N_0,
\qquad F_1= Q_0 N_1, 
\qquad F_2= Q_0 N_2,
\qquad F_3= Q_1 N.
\]
We will prove that for each $j=0,1,2,3$,
$F_j$ is a smooth map on $\Xm$, and will prove that 
for $\delta=\eph^\nu>0$ small ($\nu\in(0,1)$ fixed)
and $B_\delta$ a $\delta$-ball about $\sw$ in $\Xm$,
\begin{equation}\lbl{i.Fs}
\|F_j(\sw)\|_\Xm\le \delta, \qquad
\sup_{\us\in B_\delta}\|F_j'(\us)\|_{\LL(\Xm)}\le \delta, \qquad
j=1,2,3.
\end{equation}
(Here $\|\cdot\|_{\LL(\Xm)}$ denotes the operator norm on $\Xm$.)
The estimates in
Theorem~\ref{t.swa} follow directly from Lemma A.1 of \cite{FP1}
by these estimates and the fact that $F_0$ is smooth and 
$I-F_0'(\sw)$ has bounded inverse. 

We proceed to prove the estimates in \qref{i.Fs}.
It is clear that each 
$N_j$ ($j=0,1,2,3$) is smooth in $B_\delta$.
Also it is not hard to see 
that for some constant independent of $\eph$,
the remainder term in the nonlinearity satisfies
\begin{equation}
\|N_2(\sw)\|_\Xm + 
\sup_{\us\in B_\delta} 
\|N_2'(\us)\|_{\LL(\Xm)} \le K\eph
\end{equation}
Since $Q_0$ is bounded on $\Xm$, the estimates \qref{i.Fs} hold for
$F_2$ with $\delta=K\eph$.  We now need one more symbol estimate. 
\begin{lemma}\lbl{L.Sbd}
Let $\nu\in(0,\half)$. Then there exists $K$ such that for $\eph$
sufficiently small,
\[
\|Q_1\|_{\LL(\Xm)} = 
\sup_{k\in\R} |\hat Q_1(k+i\aaw)| \le K\eph^{2\nu}. 
\]
\end{lemma}

\noindent{\it Proof.} We fix $\nu\in(0,\half)$ and consider separately 
the low-frequency case $|\eph\xi|<3\eph^\nu$ 
and its high-frequency complement, with $\xi=k+i\aaw$.

1. Consider first the regime $|\eph\xi|<3\eph^\nu$. Let
\[
D_0 = \hat Q_0(\xi)\inv = 1+\sfrac13\xi^2, \qquad 
E(\xi)=\hat Q_0(\xi)\inv-\hat Q(\xi)\inv = \xi^2O(\eph^2\xi^2) ,
\]
so that 
\[
\hat Q_1(\xi) = \hat Q(\xi)-\hat Q_0(\xi)=
\frac1{D_0}\frac{E/D_0}{1-E/D_0}.
\]
From the Taylor expansion used in \qref{e.Sexp} we infer that
\[
|\hat Q_1(\xi)| \le K\left|\frac{E}{D_0}\right| 
\le \frac{K\eph^{2\nu}|\xi|^2}{|D_0|}\le K\eph^{2\nu}.
\]
2. In the regime $|\eph\xi|>3\eph^\nu$, we estimate $\hat Q_0$ and $\hat
Q$ separately. First, for $\eph<1$ we have $|\xi|^2>9$ and
consequently $|\hat Q_0(\xi)|$ is handled by the estimate
\begin{equation}\lbl{i.S0}
|\hat Q_0(\xi)| \le \frac{6}{|\xi|^2} \le 2\eph^{2-2\nu}.
\end{equation}
It remains to bound $|\hat Q(\xi)|$. We calculate that for 
$\xi=k+i\aaw$ with $\eph$ small,
\begin{eqnarray}
\re \frac{\tanh\eph\xi}{\eph\xi} &=& 
\frac{k\sinh 2\eph k+\aaw\sin 2\eph\aaw}{\cosh 2\eph k+\cos 2\eph\aaw}
\frac{\eph\inv}{k^2+\aaw^2}
\nonumber \\
&\le&
\frac{\sec\eph\aaw}{\eph k}
\frac{\sinh 2\eph k}{\cosh 2\eph k+1}+
\frac{2\aaw^2 }{|\xi|^2}
\nonumber\\
&=&\sec{\eph\aaw}\frac{\tanh\eph k}{\eph k}+
\frac{2\aaw^2 }{|\xi|^2}
\le (1+\eph^2\aaw^2)(1-\eph^{2\nu})+ \eph^{2-2\nu}.
\end{eqnarray}
This implies that 
\[
\re\hat Q(\xi)\inv \ge \eph^{-2}\gamma(1-
(1+\eph^2\aaw^2)(1-\eph^{2\nu})- \eph^{2-2\nu}) \ge \half
\eph^{2\nu-2}
\]
for small enough $\eph$ (since $2\nu<2-2\nu$), and hence 
\begin{equation}\lbl{i.Sbd}
|\hat Q(\xi)|\le 2\eph^{2-2\nu}.
\end{equation}
This finishes the proof of the Lemma.

From this Lemma, the estimates \qref{i.Fs} for $F_3$ clearly follow with 
$\delta= K\eph^{2\nu}$.

It remains to prove \qref{i.Fs} for $F_1$, with $\delta=K\eph^2$.
To do this it is convenient to note that the operator $Q_0$ 
gains regularity---it is a bounded map from $\Xmm$ to $\Xm$. 
Since $N_1$ is quadratic, it then suffices to prove that for
some constant $K$ independent of $\eph$, we have
\begin{equation}\lbl{i.N1}
\|
(\tanh\eph D\,\us_1)
(\tanh\eph D\,\us_2)
\|_{\Xmm} \le K\eph^2 \|\us_1\|_\Xm \|\us_2\|_\Xm
\end{equation}
for all $\us_1$, $\us_2\in \Xm$.
But this follows easily since the bilinear product map is
continuous from $\Ymm\times\Ymm$ to $\Xmm$, and 
\begin{equation}
\sup_{k\in\R} \frac{|\tanh\eph(k+i\aaw)|}{\ip{k}} \le 
\eph
\sup_{k\in\R} \frac{|\tanh\eph(k+i\aaw)|}{\ip{\eph k}} \le 
K\eph
\end{equation}
which implies that for all $\us\in\Xm$, 
\[
\|(\tanh\eph D\,\us)
\|_{\Ymm} \le K\eph \|\us\|_\Xm  .
\]
This finishes the proof of the estimates in \qref{i.Fs}.

That the fixed point is a smooth function of $\eph$ 
is a standard consequence of the easily verified fact that 
the map $(\eph,\us)\mapsto F(\us)$ is smooth.

\section{Neutral modes and adjoints}

Here we verify that the translational and wave-speed
solitary-wave degrees of freedom naturally yield two independent
elements of the generalized kernel of $\Ac$ in $L^2_a$, and 
we demonstrate that the symplectic orthogonality conditions  
\qref{c.symp} transform precisely to the condition that the initial data 
for the linearized equations
lie in the spectral complement to this generalized kernel.

1. Recall that by Theorem~\ref{t.sw}, we have a smooth family of solitary waves
$(\eta,U)$ that are solutions of the equations \qref{e.Us1}.
We exploit invariance with respect to translation by 
differentiating in $x$ to obtain an eigenfunction 
of \qref{E.evp1} corresponding to $\lam=0$. 
This is slightly tricky due to the meaning of the
variable $\phi_1$ in \qref{E.evp1}.  
We claim that the eigenfunction has the form
\begin{equation}\lbl{e.efnc1}
z_1=\pmat{\eta_1 \\ \phi_1 }=
\pmat{\eta_x \\  \phi_x},
\end{equation}
where $\phi_x=\Phi_x-v\eta_x=u$ is evaluated on the surface $(x,\eta(x))$.
To justify this statement, we note that 
\begin{equation}\lbl{e.Deta}
(1-u)\eta_x = \Heta u = -v.
\end{equation}
To see this, recall from \qref{e.Us} that $\eta=\Psi=\Heta\Phi$,
and this equation continues to hold for translated wave profiles.
Differentiating with respect to the translation parameter we have
$\dot\eta=\eta_x=\dot\Psi$ and $\dot\Phi=\Phi_x=\phi_x+v\eta_x$. 
Then \qref{e.Deta} follows from the linearization formula
\qref{E.dV}.


Using \qref{e.Deta} together with direct differentiation of \qref{e.Us1} 
(as in \qref{e.lin2}, noting $V=\eta_x$)  yields
\begin{equation}
\pmat{-\Dx(1-u) & \Dx\Heta \\ \gamma-(1-u)v' & -(1-u)\D_x}
\pmat{\eta_x \\  \phi_x} = 0.
\end{equation}
Thus $\Ac_\eta z_1=0$.
Carrying out the transformations
\qref{d.eU2}, \qref{d.eu3}, \qref{d.4}
that lead to \qref{e.evfinal}, we let
\begin{equation}
z_4 = 
\pmat{1 &-1\cr 1 & 1}
\pmat{\tig q\zst &0 \\ 0&\SC \zsh}
\pmat{\eta_x \\  \phi_x}
=\pmat{1 &-1\cr 1 & 1}
\pmat{\tig q (i\tanh D)\uu \\ \SC \zsh u}.
\end{equation}
Due to the regularity from Theorem~\ref{t.swa}
and the formulae \qref{e.uvomeg} and \qref{d.pq}, 
$z_4\in (H^1_a\cap H^1_{-a})^2$ with
\[
\Ac z_4=0.
\]

2.  Next we exploit wave-speed variation to find a
generalized eigenfunction in $L^2_a$ for $\lambda=0.$
To calculate this, it is convenient to unscale the wave speed and 
keep $\gamma$ at a fixed value $\gah$ when computing variations.
For some fixed $\ch>\sqrt{gh}$ set
\begin{equation}
\gah = \frac{gh}{\ch^2}, \quad 
\tic = \frac{c}{\hat c}=\sqrt{\frac{\gah}{\gamma}}, \quad
\etas(x;\tic)=\eta(x;\gamma),\quad
\Us(x;\tic)=\tic U(x;\gamma), \quad
\end{equation}
\begin{equation}\lbl{d.Psp}
\Phis^+(x;\tic)=\tic\Phi^+(x;\gamma) = \Dx\inv\Us = 
\tic \int_{+\infty}^{\zeta\inv(x)} \uu(s;\gamma)\,ds.
\end{equation}
Then the unscaled solitary-wave profile $(\etas,\Phis^+)$ 
is a smooth function of $\tic$ that takes values in 
$H^1_a\times H^{3/2}_a$ and satisfies
\begin{eqnarray}
-\tic\etas + \HH_{\etas} \Phis^+ = 0,
\quad
-\tic\Us+\gah\etas+\frac12(\Us,\Vs)M(\etas)\inv(\Us,\Vs)^T=0,
\end{eqnarray}
where $\Vs=\D_x \HH_{\etas}\Phis^+=\tic\D_x\etas$. 
For $\tic=1$ this simply corresponds
to \qref{e.Us1}-\qref{e.Us} with $\gamma=\gah$.

We differentiate with respect to $\tic$, then set $\tic=1$
and drop the star subscripts.
Denoting the $\tic$-derivative by the subscript $c$,
using the linearization formula \qref{E.dV} we find 
$\dc V=\Dx\eta+\Dx\dc\eta$ and
\[
-\eta-\dc\eta+\Heta(\dc\Phi^+-v\dc\eta)+u\dc\eta
=0,
\]
\[
-U-\dc U+\gamma\dc\eta+u\dc U+v(\Dx\eta+\Dx\dc\eta)
-uv\Dx\dc\eta = 0.
\]
Since $U-\eta_xv=u=\phi_x$, this yields
\begin{equation}\lbl{e.pcp}
\pmat{-\Dx(1-u) & \Dx\Heta \\ \gamma-(1-u)v' & -(1-u)\D_x}
\pmat{\dc \eta\\ \dc\phi^+} = 
\pmat{\eta_x \\ \phi_x},
\qquad \dc\phi^+=\dc\Phi^+-v\dc\eta.
\end{equation}
With
\begin{equation}\lbl{d.z1}
y_1 = \pmat{\dc \eta\\ \dc\phi^+}  ,
\qquad y_4 = 
\pmat{1 &-1\cr 1 & 1}
\pmat{\tig q\zst & 0 \\ 0& \SC\zsh } y_1,
\end{equation}
we have $-\Ac_\eta y_1=z_1$, and find that 
$y_4\in (H^1_a)^2$  with
\begin{equation}
-\Ac y_4 = z_4.
\end{equation}

{\bf Adjoint modes.} 
It is a standard fact of operator theory that the space $\bar Y_a$,
the kernel of the spectral projection $P_0$ in \qref{d.P0},
is the subspace annihilated by the generalized kernel 
of the adjoint $\Ac^*$.  This generalized kernel is two-dimensional
(since the generalized kernel of $\Ac$ is), 
and we aim to show that the annihilation conditions correspond 
to the symplectic orthogonality conditions \qref{c.symp}.

We will work with the Banach space dual $L^2_{-a}$ of $L^2_a$,
and note that for the Fourier multiplier 
$\SC=\sqrt{-\tig D\tanh D}$, the adjoint is given formally by $S^*=S$ 
acting in $L^2_{-a}$. To see this, take smooth test functions $f$ and $g$
and write $f_a=\eax f$, $g_{-a}=\emax g$, and $\SC_a=\eax\SC\emax$.
Then since the symbol of $\SC$ satisfies
$\overline{\SC(k+ia)}=\SC(k-ia)$ we have
\begin{eqnarray*}
\intR (\SC f)\overline{g}\,dx 
&=& \intR \SC_a f_a \overline{g_{-a}}\,dx
= \intR \SC(k+ia)\hat f_a(k)\overline{\hat g_{-a}(k)}\,\frac{dk}{2\pi}
\\ &=& \intR f_a \overline {\SC_{-a} g_{-a}}\,dx
= \intR f \overline{\SC g}\,dx.
\end{eqnarray*}

To describe the generalized kernel of the adjoint $\Ac^*$,
first note that with the definition
\begin{equation}\lbl{e.Pcm}
\Phis^-(x;\tic) = \tic \Phi^-(x;\gamma)=
\tic\int_{-\infty}^{\zeta\inv(x)}\uu(s;\gamma)\,ds,
\end{equation}
we can repeat the arguments leading up to \qref{e.pcp} 
with $\Phi^-$ and $H^s_{-a}$ replacing $\Phi^+$ and $H^s_a$.
Then \qref{e.pcp} holds with $\dc\phi^-=\dc\Phi^- -v\dc\eta$ 
replacing $\dc\phi^+$.

Next, it is convenient to work with the variables used in 
\qref{d.eu3}-\qref{e.ev3} and note that 
\[
\Ac_*= 
\pmat{1 &-1\cr 1 & 1} 
\Ac_3
\pmat{1 &1\cr -1 & 1}
\frac12,
\qquad
\Ac_3 = 
\pmat{ q\D pq\inv & -q\SC \cr -\SC q &\SC p\D\SC\inv}.
\]
The operator $\Ac_3$ acts in $(L^2_a)^2$.
As can be found by transformation from the original 
canonical Hamiltonian structure, this operator admits the 
factorization 
\begin{equation}\lbl{e.JL3}
\Ac_3=\JJ\LL, \qquad 
\JJ = 
\pmat{0 & q\SC\cr -\SC q & 0}, 
\quad
\LL = \pmat{1 & -q\inv p\D\SC\inv \cr \SC\inv\D p q\inv & -1}.
\end{equation}
The adjoints are given by $\JJ^*=-\JJ$, $\LL^*=\LL$, and
$\Ac_3^*=-\LL\JJ$, acting in $(L^2_{-a})^2$.
Then with
\begin{equation}\lbl{d.z3}
z_3 = 
\pmat{\tig q\zst & 0 \\ 0& \SC\zsh } 
\pmat{\eta_x\\ \phi_x},
\qquad y_3=
\pmat{\tig q\zst & 0 \\ 0& \SC\zsh } 
\pmat{\dc\eta \\ \dc\phi^-} ,
\end{equation}
\begin{equation}
z_3^* = \tig\inv\JJ\inv z_3= 
\pmat{ -(\tig q)\inv \zsh\phi_x \\
\SC\inv \zst\eta_x},
\qquad y_3^* = \tig\inv\JJ\inv y_3
= \pmat{ -(\tig q)\inv \zsh\phi^-_c \\
\SC\inv \zst\eta_c},
\end{equation}
one can check directly that
$z_3^*$, $y_3^*\in (H^1_{-a})^2$ and 
\begin{eqnarray*}
\Ac_3^* z_3^* &=& 
\pmat{-q\inv p\D q & -q\SC \cr -\SC q & -\SC\inv\D p\SC}
\pmat{ -(\tig q)\inv \zsh\phi_x \\
\SC\inv \zst\eta_x}
\\
&=& \pmat{ (\tig q)\inv &0 \cr 0 &\SC\inv}
\pmat{ -p\D & \tig q^2\cr D\tanh D & \D p}
\pmat{\zsh \phi_x\cr \zst \eta_x} = 0,
\end{eqnarray*}
and similarly $-\Ac_3^*y_3^* = z_3^*$.
Thus the generalized kernel of $\Ac_3^*$ is the span of
$z_3^*$ and $y_3^*$.

Now, corresponding to an arbitrary element
$\dot z_1=(\dot \eta,\dot\phi)\in Z_a=L^2_a\times
H^{1/2}_a$ is
$\dot z_3=(\tig q\zst\dot\eta,\SC\zsh\dot\phi)\in (L^2_a)^2$.
Then the conditions that $z_3^*$ and $y_3^*$ annihilate $\dot z_3$
transform as follows:
\begin{equation}
0= -\ip{\dot z_3,z_3^*} =
\intR (\tig q\zst \dot\eta)
\overline{(\tig q)\inv \zsh\phi_x}-
(S\zsh\dot\phi)\overline{\SC\inv\zst\eta_x} \,d\xf
=
\intR \dot\eta \phi_x - \dot\phi \eta_x \, dx,
\end{equation}
\begin{equation}
0= -\ip{\dot z_3,y_3^*} =
\intR (\tig q\zst \dot\eta)
\overline{(\tig q)\inv \zsh\phi^-_c}-
(S\zsh\dot\phi)\overline{\SC\inv\zst\eta_c} \,d\xf
=
\intR \dot\eta \phi^-_c - \dot\phi \eta_c \, dx.
\end{equation}
This shows that the symplectic orthogonality conditions
\qref{c.symp} transform to the precise condition 
that the initial data for the semigroup $e^{\Ac t}$ lie
in the space $\bar Y_a=\ker P_0$ that is the spectral complement 
of the generalized kernel of $\Ac$.

\section{Characteristic values for the KdV bundle}

Here we provide the proof of Proposition~\ref{p.kdv} concerning 
the characteristic values of the bundle 
\[
W_0(\lam)= I + (\lam-\half\D+\sfrac16\D^3)\inv\D(\sfrac32\sw),
\qquad \sw = \sech^2(\sqrt3 x/2).
\]
(In this section we will drop the tilde on $\lam$ for convenience.)
For $\re\lam>\beta=\half\aaw(1-\frac13\aaw^2)$, 
the weight-transformed operator 
\[
\eapx (I-W_0(\lam))\emapx=
(-\lam+\half(\D-\aaw)+\sfrac16(\D-\aaw)^3)\inv(\D-\aaw)
(\sfrac32\sw )
\]
is compact on $L^2$ and has range in $H^1$, 
due to estimates similar to \qref{i.mo1}-\qref{i.mobd1}. 
On $L^2_\aaw$, therefore, $W_0(\lam)$ is Fredholm of index zero. 
If $W_0(\lam)f=0$ for some nonzero $f\in L^2_\aaw$, 
then $\eapx f\in H^m$ for all $m$ by an easy bootstrapping argument, 
and $f$ satisfies the ordinary differential equation
\begin{equation}\lbl{e.weq}
(\lam -\sfrac12\D+\sfrac16\D^3)f+\D\left(\sfrac32\sw f\right) = 0.
\end{equation}
By standard results, such an equation has a solution 
$f\sim e^{\mu x}$ as $x\to\infty$ for each $\mu$ that satisfies
\begin{equation}\lbl{e.murt}
\lam-\half\mu+\sfrac16\mu^3=0.
\end{equation}
For $\lam>0$ large, this equation has 
one root with $\re\mu<-\aaw$  and two with $\mu$ with $\re\mu>-\aaw$.
With $\mu=-\aaw+it$ ($t\in\R$), the curve 
\[
t\mapsto \mu-\sfrac13\mu^3 = 
-\aaw+\sfrac13\aaw^3-\aaw t^2 + i(t+\sfrac13t^3-\aaw^2 t)
\]
has increasing imaginary part and real part less than 
$-\aaw+\sfrac13\aaw^3=-2\beta$.
For $\re\lam>-\beta$, then, \qref{e.murt} has a unique and simple root 
satisfying $\re\mu<-\aaw$, hence \qref{e.weq} has a unique solution  
(up to a constant factor)  satisfying $\eapx f\to 0$ as $x\to\infty$.
In particular, one may check explicitly (and easily by computer) that 
\begin{equation}\lbl{d.f}
f = \Dx \left( e^{\mu x}
((\sqrt3+\mu)^2-
(\sqrt3+\mu+\mu e^{\sqrt3x})
\sqrt3\sech^2(\sqrt3 x/2)
)
\right).
\end{equation}
Since $\re\mu<-\aaw$, clearly $\eapx f\in L^2$ is impossible unless 
$\sqrt3+\mu=0$, meaning $\lam=0$ and $f=\Dx\sw$.

From the analysis so far, we see that the kernel of $W_0(0)$ in $L^2_\aaw$ 
is one-dimensional. To finish the proof, we need to show that there is 
a root vector at 0 with order 2, and no root vector of order 3.
Any root vector $f(\lam)$ at $\lam=0$ may be taken in the form
$f(\lam)=f_0+\lam f_1+\lam^2 f_2+O(\lam^3)$ with $f_0=\Dx\sw$.
And $W_0(\lam)=W_0+ W_0'\lam+\half W_0''\lam^2+O(\lam^3)$ where
$W_0=W_0(0)$ and
\[
W_0'=
-(-\half\D+\sfrac16\D^3)^{-2}\D(\sfrac32\sw),
\qquad
\half W_0''=(-\half\D+\sfrac16\D^3)^{-3}\D(\sfrac32\sw).
\]
To find a root vector of order 2, it suffices to find 
$f_1\in L^2_\aaw$ such that $W_0f_1+W_0'f_0=0$. Since $W_0f_0=0$ we have 
\[
W_0'f_0 = (-\half\D+\sfrac16\D^3)\inv f_0,
\]
so it suffices to find $f_1$ such that
\[
(-\half\D+\sfrac16\D^3)f_1+ \D(\sfrac32\sw f_1) + f_0 = 0.
\]
Such a function can be found by differentiating the equation
satisfied by the KdV wave profile with respect to wave speed.
The function 
$ \varphi_b(x)=b\sech^2 \sqrt{3b}x/2$
satisfies $\varphi_1=\sw$ and
\[
(-\half b \Dx +\sfrac16 \Dx^3)
\varphi_b + \Dx(\sfrac34 \varphi_b^2)=0.
\]
Differentiating with respect to $b$ and setting $b=1$, we find that 
\[
(-\half \Dx +\sfrac16 \Dx^3+\sfrac32\Dx\sw)\D_b\varphi_1
-\sfrac12\Dx\varphi_1=0.
\]
From this we see that the choice $f_1=-\half\D_b\varphi_1$ works and 
yields a root vector of order 2.  This choice is unique up to adding
a scalar multiple of $f_0$.

To show that there is no root vector of order greater than 2, it suffices
to show that with $f_1$ as above, there is no $f_2\in L^2_\aaw$ such
that 
\begin{equation}\lbl{e.w2}
W_0f_2+W_0'f_1+\half W_0'' f_0 = 0.
\end{equation}
If such an $f_2$ exists, then a bootstrapping argument involving the 
decay estimate \qref{i.mobd1} shows that $e^{\aaw x}f_2\in H^m$ for all $m$.
Because of the equations satisfied by $f_1$ and $f_0$, we find that 
\[
W_0'f_1+\half W_0''f_0 =  
(-\half \Dx +\sfrac16 \Dx^3)\inv f_1.
\]
Therefore $f_2$ must be a smooth solution of 
\begin{equation}\lbl{e.f2e}
(-\half \Dx +\sfrac16 \Dx^3+\sfrac32\Dx\sw)f_2 + f_1 = 0.
\end{equation}
Now, the function $\varphi_1=\sw$ has $\emapx \varphi_1$ in $H^m$ for all $m$ and 
satisfies 
\[
(-\half \Dx +\sfrac16 \Dx^3+\sfrac32\sw\Dx)\varphi_1 = 0.
\]
Multiplying \qref{e.f2e} by $\varphi_1$ and integrating by parts, we find that the 
terms involving $f_2$ vanish. Thus, for $f_2$ to exist, it is necessary that
$\intR \varphi_1\D_b\varphi_1 = 0$. But
\[
\intR \varphi_1\D_b\varphi_1\,dx = \frac{d}{db}\intR \frac12\varphi_b^2\,dx =
\frac{d}{db} b^{3/2}\intR\varphi_1^2\,dx>0.
\]
Hence, $f_2$ cannot exist as required, and this proves that the characteristic
value $\lam=0$ has null multiplicity 2.

\section*{Acknowledgments}
RLP is grateful for discussions regarding this work with
Michael Weinstein, and for partial support during the completion
of this work by the Oxford Centre for Nonlinear PDE.
This material is based upon work supported by the National Science
Foundation under grant nos.\ DMS 06-04420, 09-05723, and 08-07597,
and by the Center for Nonlinear Analysis under 
NSF grant DMS 06-35983.  

\hide{\singlespace}
\baselineskip=12pt

\bibliographystyle{siam}
\bibliography{wave}

\begin{thebibliography}{10}

\bibitem{Be77}
{\sc J.~T. Beale}, {\em The existence of solitary water waves}, Comm. Pure
  Appl. Math., 30 (1977), pp.~373--389.

\bibitem{Be91}
\leavevmode\vrule height 2pt depth -1.6pt width 23pt, {\em Exact solitary water
  waves with capillary ripples at infinity}, Comm. Pure Appl. Math., 44 (1991),
  pp.~211--257.

\bibitem{Ben72}
{\sc T.~B. Benjamin}, {\em The stability of solitary waves}, Proc. Roy. Soc.
  (London) Ser. A, 328 (1972), pp.~153--183.

\bibitem{BeLf}
{\sc J.~Bergh and J.~L{\"o}fstr{\"o}m}, {\em Interpolation spaces. {A}n
  introduction}, Springer-Verlag, Berlin, 1976.
\newblock Grundlehren der Mathematischen Wissenschaften, No. 223.

\bibitem{Bo75}
{\sc J.~Bona}, {\em On the stability theory of solitary waves}, Proc. Roy. Soc.
  London Ser. A, 344 (1975), pp.~363--374.

\bibitem{BoSa89}
{\sc J.~L. Bona and R.~L. Sachs}, {\em The existence of internal solitary waves
  in a two-fluid system near the {K}d{V} limit}, Geophys. Astrophys. Fluid
  Dynam., 48 (1989), pp.~25--51.

\bibitem{Bou1871a}
{\sc J.~Boussinesq}, {\em Th\'eorie de l'intumescence liquide appel\'ee onde
  solitaire ou de translation se propageant dans un canal rectangulair}, C. R.
  Acad. Sci. Paris, 72 (1871), pp.~755--759.

\bibitem{Bou1871b}
\leavevmode\vrule height 2pt depth -1.6pt width 23pt, {\em Th\'eorie
  g\'en\'erale des mouvements qui sont propag\'es dans un canal rectangulaire
  horizontal}, C. R. Acad. Sci. Paris, 73 (1871), pp.~256--260.

\bibitem{Bou1872}
\leavevmode\vrule height 2pt depth -1.6pt width 23pt, {\em Th\'eorie des ondes
  et des remous qui se propagent le long d'un canal rectangulaire horizontal,
  en communiquant au liquide contenu dans ce canal des vitesses sensiblement
  pareilles de la surface au fond}, J. Math. Pures et Appliq., 17 (1872),
  pp.~55--108.

\bibitem{Bou1877}
\leavevmode\vrule height 2pt depth -1.6pt width 23pt, {\em Essai sur la
  th\'eorie des eaux courantes}, M\'emoires pr\'esent\'es par divers savants
  \'a l'Acad\'emie des Sciences Inst. France (s\'eries 2), 23 (1877),
  pp.~1--680.

\bibitem{Bu04}
{\sc B.~Buffoni}, {\em Existence and conditional energetic stability of
  capillary-gravity solitary water waves by minimisation}, Arch. Ration. Mech.
  Anal., 173 (2004), pp.~25--68.

\bibitem{Craig}
{\sc W.~Craig}, {\em An existence theory for water waves and the {B}oussinesq
  and {K}orteweg-de {V}ries scaling limits}, Comm. Partial Differential
  Equations, 10 (1985), pp.~787--1003.

\bibitem{FrHy54}
{\sc K.~O. Friedrichs and D.~H. Hyers}, {\em The existence of solitary waves},
  Comm. Pure Appl. Math., 7 (1954), pp.~517--550.

\bibitem{FP1}
{\sc G.~Friesecke and R.~L. Pego}, {\em Solitary waves on {FPU} lattices. {I}.
  {Q}ualitative properties, renormalization and continuum limit}, Nonlinearity,
  12 (1999), pp.~1601--1627.

\bibitem{FP2}
\leavevmode\vrule height 2pt depth -1.6pt width 23pt, {\em Solitary waves on
  {FPU} lattices. {II}. {L}inear implies nonlinear stability}, Nonlinearity, 15
  (2002), pp.~1343--1359.

\bibitem{FP3}
\leavevmode\vrule height 2pt depth -1.6pt width 23pt, {\em Solitary waves on
  {F}ermi-{P}asta-{U}lam lattices. {III}. {H}owland-type {F}loquet theory},
  Nonlinearity, 17 (2004), pp.~207--227.

\bibitem{FP4}
\leavevmode\vrule height 2pt depth -1.6pt width 23pt, {\em Solitary waves on
  {F}ermi-{P}asta-{U}lam lattices. {IV}. {P}roof of stability at low energy},
  Nonlinearity, 17 (2004), pp.~229--251.

\bibitem{GGK90}
{\sc I.~Gohberg, S.~Goldberg, and M.~A. Kaashoek}, {\em Classes of linear
  operators. {V}ol. {I}}, vol.~49 of Operator Theory: Advances and
  Applications, Birkh\"auser Verlag, Basel, 1990.

\bibitem{GS71}
{\sc I.~C. Gohberg and E.~I. Sigal}, {\em An operator generalization of the
  logarithmic residue theorem and {R}ouch\'e's theorem}, Mat. Sb. (N.S.),
  84(126) (1971), pp.~607--629.

\bibitem{GSS1}
{\sc M.~Grillakis, J.~Shatah, and W.~Strauss}, {\em Stability theory of
  solitary waves in the presence of symmetry. {I}}, J. Funct. Anal., 74 (1987),
  pp.~160--197.

\bibitem{GSS2}
\leavevmode\vrule height 2pt depth -1.6pt width 23pt, {\em Stability theory of
  solitary waves in the presence of symmetry. {II}}, J. Funct. Anal., 94
  (1990), pp.~308--348.

\bibitem{GW10}
{\sc M.~D. Groves and E.~Wahl{\'e}n}, {\em On the existence and conditional
  energetic stability of solitary water waves with weak surface tension}, C. R.
  Math. Acad. Sci. Paris, 348 (2010), pp.~397--402.

\bibitem{HS02a}
{\sc M.~Haragus and A.~Scheel}, {\em Finite-wavelength stability of
  capillary-gravity solitary waves}, Comm. Math. Phys., 225 (2002),
  pp.~487--521.

\bibitem{Kato}
{\sc T.~Kato}, {\em Perturbation theory for linear operators}, Die Grundlehren
  der mathematischen Wissenschaften, Band 132, Springer-Verlag New York, Inc.,
  New York, 1966.

\bibitem{KdV1895}
{\sc D.~J. Korteweg and G.~de~Vries}, {\em On the change of form of long waves
  advancing in a rectangular canal, and on a new type of long stationary
  waves}, Phil. Mag. (Ser. 5), 39 (1895), pp.~422--443.

\bibitem{La05}
{\sc D.~Lannes}, {\em Well-posedness of the water-waves equations}, J. Amer.
  Math. Soc., 18 (2005), pp.~605--654 (electronic).

\bibitem{Lav54}
{\sc M.~A. Lavrent'ev}, {\em On the theory of long waves; {A} contribution to
  the theory of long waves}, Amer. Math. Soc. Transl., 102 (1954), pp.~3--50.

\bibitem{Lin09}
{\sc Z.~Lin}, {\em On linear instability of 2{D} solitary water waves}, Int.
  Math. Res. Not. IMRN,  (2009), pp.~1247--1303.

\bibitem{Mi01}
{\sc A.~Mielke}, {\em On the energetic stability of solitary water waves}, R.
  Soc. Lond. Philos. Trans. Ser. A Math. Phys. Eng. Sci., 360 (2002),
  pp.~2337--2358.
\newblock Recent developments in the mathematical theory of water waves
  (Oberwolfach, 2001).

\bibitem{Miles}
{\sc J.~W. Miles}, {\em The {K}orteweg-de {V}ries equation: A historical
  essay}, J. Fluid Mech., 106 (1981), pp.~131--147.

\bibitem{MiWe96}
{\sc J.~R. Miller and M.~I. Weinstein}, {\em Asymptotic stability of solitary
  waves for the regularized long-wave equation}, Comm. Pure Appl. Math., 49
  (1996), pp.~399--441.

\bibitem{MZ09}
{\sc T.~Mizumachi}, {\em Asymptotic stability of lattice solitons in the energy
  space}, Comm. Math. Phys., 288 (2009), pp.~125--144.

\bibitem{MP08}
{\sc T.~Mizumachi and R.~L. Pego}, {\em Asymptotic stability of {T}oda lattice
  solitons}, Nonlinearity, 21 (2008), pp.~2099--2111.

\bibitem{Pe85}
{\sc R.~L. Pego}, {\em Compactness in {$L^2$} and the {F}ourier transform},
  Proc. Amer. Math. Soc., 95 (1985), pp.~252--254.

\bibitem{PSun}
{\sc R.~L. Pego and S.~M. Sun}, {\em On the transverse linear instability of
  solitary water waves with large surface tension}, Proc. Roy. Soc. Edinburgh
  Sect. A, 134 (2004), pp.~733--752.

\bibitem{PW94}
{\sc R.~L. Pego and M.~I. Weinstein}, {\em Asymptotic stability of solitary
  waves}, Comm. Math. Phys., 164 (1994), pp.~305--349.

\bibitem{Pr85}
{\sc J.~Pr{\"u}ss}, {\em On the spectrum of {$C_{0}$}-semigroups}, Trans. Amer.
  Math. Soc., 284 (1984), pp.~847--857.

\bibitem{Ray1876}
{\sc L.~Rayleigh}, {\em On waves}, Phil. Mag. (Ser. 5), 1 (1876), pp.~257--279.

\bibitem{RT}
{\sc F.~Rousset and N.~Tzvetkov}, {\em Transverse instability of the line
  solitary water waves}.
\newblock Preprint at {\tt arXiv:0906.2487v2}.

\bibitem{Russell}
{\sc J.~S. Russell}, {\em Report on waves}, in Report of the 14th Meeting of
  the British Association for the Advancement of Science, London, U.K., 1844,
  John Murray, pp.~311 -- 390.

\bibitem{ScWa00b}
{\sc G.~Schneider and C.~E. Wayne}, {\em The long-wave limit for the water wave
  problem. {I}. {T}he case of zero surface tension}, Comm. Pure Appl. Math., 53
  (2000), pp.~1475--1535.

\bibitem{Sun99}
{\sc S.~M. Sun}, {\em Non-existence of truly solitary waves in water with small
  surface tension}, R. Soc. Lond. Proc. Ser. A Math. Phys. Eng. Sci., 455
  (1999), pp.~2191--2228.

\bibitem{Wu2d}
{\sc S.~Wu}, {\em Well-posedness in {S}obolev spaces of the full water wave
  problem in {$2$}-{D}}, Invent. Math., 130 (1997), pp.~39--72.

\bibitem{Wu09}
\leavevmode\vrule height 2pt depth -1.6pt width 23pt, {\em Almost global
  wellposedness of the 2-{D} full water wave problem}, Invent. Math., 177
  (2009), pp.~45--135.

\bibitem{Za68}
{\sc V.~E. Zakharov}, {\em Stability of periodic waves of finite amplitude on
  the surface of a deep fluid}, Zh. Prikl. Mekh. Fiz., 9 (1968), pp.~86--94.

\end{thebibliography}

\end{document}